\documentclass{article}

\title{Frobenius algebra objects in\\Temperley-Lieb categories at roots of unity}
\author{Joseph Grant\thanks{joseph.grant@gmail.com} {}
 and Mathew Pugh\thanks{School of Mathematics, Cardiff University, Abacws, Senghennydd Road, Cardiff, CF24 4AG}}
\date{\vspace{-2em}}

\usepackage[UKenglish]{isodate}
\cleanlookdateon

\usepackage{sfmath}
\usepackage{tikz}

\usepackage{latexsym}
\usepackage{amsmath}
\usepackage{amssymb}
\usepackage{amsthm}
\usepackage[shortlabels]{enumitem}

\usepackage{ dsfont } 
\usepackage{graphicx}

\usepackage{hyperref}  

\input xy
\xyoption{all}

\newcommand{\ti}{\mathds{1}}
\newcommand{\id}{1} 

\newcommand{\R}{\mathbb{R}}
\newcommand{\Cnum}{\mathbb{C}}
\newcommand{\N}{\mathbb{N}}
\newcommand{\tl}{\widetilde{TL}}
\newcommand{\tlk}{\tl_k}

\newcommand{\then}{\Rightarrow}

\newcommand{\Ccat}{\mathcal{C}}
\newcommand{\D}{\mathcal{D}}
\newcommand{\M}{\mathcal{M}}
\newcommand{\kk}{\mathds{k}}
\newcommand{\Vect}{{\operatorname{Vec}\nolimits}}
\newcommand{\ob}{\operatorname{ob}\nolimits}

\newcommand{\End}{\operatorname{End}\nolimits}
\newcommand{\spn}{\operatorname{span}\nolimits}

\newcommand{\into}{\hookrightarrow}

\newcommand{\Hom}{\operatorname{Hom}\nolimits}
\newcommand{\Fun}{\operatorname{Fun}\nolimits}
\newcommand{\mMod}{\operatorname{\!-mod}}
\newcommand{\Modm}{\operatorname{mod-}\!}
\newcommand{\mModm}{\!\operatorname{-mod-}\!}

\newcommand{\da}{\text{-}}
\newcommand{\arr}[1]{\stackrel{#1}{\to}}
\newcommand{\larr}[1]{\stackrel{#1}{\longrightarrow}}
\newcommand{\arrr}[1]{\larr{#1}}

\newtheorem{theorem}{Theorem}[section]

\newtheorem{corollary}[theorem]{Corollary}
\newtheorem{lemma}[theorem]{Lemma}
\newtheorem{proposition}[theorem]{Proposition}
\newtheorem{definition-proposition}[theorem]{Definition-Proposition}

\newtheorem*{thm0}{Theorem}

\theoremstyle{definition}
\newtheorem{definition}[theorem]{Definition}
\newtheorem{remark}[theorem]{Remark}
\newtheorem{example}[theorem]{Example}

\begin{document}

\maketitle

\begin{abstract}
We give a new definition of a Frobenius structure on an algebra object in a monoidal category, generalising Frobenius algebras in the category of vector spaces.  Our definition allows Frobenius forms valued in objects other than the unit object, and can be seen as a categorical version of Frobenius extensions of the second kind.  When the monoidal category is pivotal we define a Nakayama morphism for the Frobenius structure and explain what it means for this morphism to have finite order.

Our main example is a well-studied algebra object in the (additive and idempotent completion of the) Temperley-Lieb category at a root of unity.  We show that this algebra has a Frobenius structure and that its Nakayama morphism has order 2.  As a consequence, we obtain information about Nakayama morphisms of preprojective algebras of Dynkin type, considered as algebras over the semisimple algebras on their vertices.

\end{abstract}

\tableofcontents

\setlength{\parindent}{0pt}
\setlength{\parskip}{1em plus 0.5ex minus 0.2ex}


\section{Introduction}

Frobenius algebras are a central algebraic structure in modern mathematics.  They first appeared in representation theory, connected to the study of group representations \cite{bn,nes}.  They have since appeared in many other areas of mathematics and physics including the study of topological quantum field theories \cite{kock}.

A Frobenius algebra comes with a distinguished automorphism called the Nakayama automorphism.  There has been work studying when these automorphisms have finite order, especially for preprojective algebras of Dynkin quivers and related algebras \cite{bbk, hi-frac, ep-nak, g-ser}.

One can obtain the Dynkin preprojective algebras $\Pi$ in the following way.  First, one constructs the Temperley-Lieb category $TL$ as the additive and idempotent completion of the $\Cnum$-linear category of diagrams on non-crossing strings between dots, modulo a circle relation with quantum parameter a fixed root of unity.  The category $TL$ has a monoidal structure given by juxtaposing diagrams and is known to be a fusion category with simple objects given by Jones-Wenzl projections.  Quotienting $TL$ by an ideal of negligible morphisms gives a monoidal category $\tl$.
Now one can construct an algebra object $\Sigma$ in $\tl$ from the Jones-Wenzl projectors and study monoidal functors from $\tl$ to categories $S\mModm S$ of bimodules over a semisimple algebra.  The image of $\Sigma$ under such a functor becomes an algebra object $\Pi$ in $S\mModm S$ whose underlying $\Cnum$-algebra is a Dynkin preprojective algebra.

The axioms for a Frobenius algebra are easily imitated in the setting of monoidal categories, so one might hope that the algebra objects $\Pi\in S\mModm S$ and $\Sigma\in\tl$ are Frobenius, but this is not the case: see Example \ref{ex:toyPi} and Lemma \ref{lem:no-classical} below.

In this article we introduce a generalised definition of a Frobenius algebra object $A$ in a rigid monoidal category $\Ccat$ which involves a ``twist'' by another object $W\in\Ccat$: we allow Frobenius forms $A\to W$ instead of just maps $A\to\ti$ to the unit object.  In Proposition \ref{prop:2ndkind} we see that these generalised Frobenius structures are related to the classical notion of Frobenius extensions of the second kind.

When our monoidal category is pivotal, i.e., it comes with natural isomorphisms from an object to its double dual, we can define a Nakayama morphism of a Frobenius algebra object: this is no longer an automorphism, but it plays a similar role in the theory.  When the object $W\in\Ccat$ has finite order, i.e., there is an isomorphism $W\otimes\cdots\otimes W\cong \ti$, it makes sense to ask about the order of the Nakayama morphism.

We show that the algebra object $\Sigma\in\tl$ satisfies our definition of a Frobenius algebra object, where the twisting object $W$ is given by the highest Jones-Wenzl projector.  Then, by using the braiding on $\tl$, we are able to determine the order of the Nakayama automorphism of $\Sigma$.
\begin{thm0}[{Lemma \ref{lem:sigma-frob} and Theorem \ref{thm:alpha^2}}]
The algebra object $\Sigma\in\tl$ is Frobenius with Nakayama morphism of order $2$.
\end{thm0}

As a corollary, we show that the algebra object $\Pi\in S\mModm S$ is also Frobenius with Nakayama morphism of order $2$.  Our approach involves a single calculation in $\tl$.  This is in contrast to the classical proof that the Nakayama automorphism of the $\Cnum$-algebra $\Pi$ has a Nakayama automorphism of order $2$, which involves separate calculations in types $A$, $D$, and $E$ \cite{bbk}.

\emph{Acknowledgements:} 
Thanks to Alastair King for explaining results in \cite{cooper}.


\section{Frobenius algebras}

\subsection{Frobenius structures}
We recall some standard notions, mainly to fix notation.  The reader can find more details about monoidal categories, algebra objects, and dualizable objects in sections 2.1, 7.8, and 2.10 of \cite{egno}.

Let $(\Ccat, \otimes, \ti)$ be a monoidal category.  Recall that an \emph{algebra} (or \emph{algebra object}) in $\Ccat$ is a triple $(A,m,u)$ where $A\in\Ccat$ is an object and $m$ and $u$ are maps $m:A\otimes A\to A$ (\emph{multiplication}) and $u:\ti\to A$ (\emph{unit}) in $\Ccat$ which satisfy the associativity and unitality conditions. 
We sometimes draw maps using string diagrams which go up the page, as follows:
\[
\begin{tikzpicture}
[auto, block/.style ={rectangle, draw=black, thick, fill=blue!10, text width=2em, align=center, rounded corners, minimum height=2em},
 sdot/.style = {circle,draw=black, minimum width=5pt, fill=black, inner sep=0pt]} ]
\node (b1) at (0,0) {$A$};
\node (b2) at (1.4,0) {$A$};
\node (u3) at (4,1.4) {$A$};
\node (t1) at (0.7,1.4) {$A$};
\node[sdot] (m) at (0.7,0.7) {};
\node[sdot] (u) at (4,0.7) {};
\node[right] at (m) {$\;m$};
\node[right] at (u) {$\;u$};
\draw[thick] (b1) -- (m);
\draw[thick] (b2) -- (m);
\draw[thick] (t1) -- (m);
\draw[thick] (u3) -- (u);
\end{tikzpicture}
\]
The unit object is denoted by the empty diagram.

An object $X\in\Ccat$ is said to be \emph{left dualizable} if it has a left dual $X^\vee$, so there exist maps ${e_X:X^\vee\otimes X\to \ti}$ (\emph{evaluation}) and $c_X:\ti\to X\otimes X^\vee$ (\emph{coevaluation}) satisfying the triangle identities.
\[
\begin{tikzpicture}
[auto, block/.style ={rectangle, draw=black, thick, fill=blue!10, text width=2em, align=center, rounded corners, minimum height=2em},
 sdot/.style = {circle,draw=black, minimum width=5pt, fill=black, inner sep=0pt]} ]
\node (b1) at (0,0) {$X^\vee$};
\node (b2) at (1,0) {$X$};
\node (t1) at (2,1.4) {$X$};
\node (t2) at (3,1.4) {$X^\vee$};
\node (a1) at (4.8,1.4) {$X$};
\node (a2) at (5.8,0.7) {$X^\vee$};
\node (a3) at (6.8,0) {$X$};
\node (t3) at (8,1.4) {$X$};
\node (b3) at (8,0) {$X$};
\node at (7.4,0.7) {$=$};
\coordinate (y1) at (4.8,0.4) {};
\coordinate (y2) at (6.8,0.9) {};
\node (c1) at (9.8,0) {$X^\vee$};
\node (c2) at (10.8,0.7) {$X$};
\node (c3) at (11.8,1.4) {$X^\vee$};
\node (t4) at (13,1.4) {$X^\vee$};
\node (b4) at (13,0) {$X^\vee$};
\node at (12.4,0.7) {$=$};
\coordinate (z1) at (9.8,1) {};
\coordinate (z2) at (11.8,0.5) {};
\draw[thick] (b1) to[out=90,in=90] (b2);
\draw[thick] (t1) to[out=-90,in=-90] (t2);
\draw[thick] (a1) to[out=-90,in=90] (y1) to[out=-90,in=-90] (a2);
\draw[thick] (a2)  to[out=90,in=90] (y2) to (a3);
\draw[thick] (b3) to (t3);
\draw[thick] (c1) to[out=90,in=-90] (z1) to[out=90,in=90] (c2);
\draw[thick] (c2)  to[out=-90,in=-90] (z2) to (c3);
\draw[thick] (b4) to (t4);
\end{tikzpicture}
\]
Similarly, $X\in\Ccat$ is \emph{right dualizable} if it has a right dual ${}^\vee X$, so there exist maps $e^r_X:X\otimes {}^\vee X\to \ti$ and $c^r_X:\ti\to {}^\vee X\otimes X$ satisfying the triangle identities.  If $X$ is left dualizable then ${}^\vee (X^\vee)\cong X$.

\begin{definition}\label{def:frob}
Suppose $A$ is left dualizable.  A \emph{Frobenius structure} on the algebra $(A,u,m)$ is a pair $(W,n)$ where $W\in\Ccat$ and $n:A\to W$ is a map in $\Ccat$ 
such that the map $\kappa=n\circ m:A\otimes A\to W$ is non-degenerate, i.e., the map
$ \varphi:A \larr{\id_A\otimes c_A} A\otimes A\otimes A^\vee\larr{\kappa\otimes \id_{A^\vee}} W\otimes A^\vee $
\[
\begin{tikzpicture}
[auto, block/.style ={rectangle, draw=black, thick, fill=blue!10, text width=2em, align=center, rounded corners, minimum height=2em},
 sdot/.style = {circle,draw=black, minimum width=5pt, fill=black, inner sep=0pt]} ]
\node (a) at (3,1.2) {$A$};
\node (b) at (2.4,2.8) {$W$};
\node (c) at (3.6,2.8) {$A^\vee$};
\node[sdot] (phi) at (3,2) {};
\node[right] at (phi) {$\;\varphi$};
\node at (4.2,2) {$=$};
\node (d) at (4.5,0) {$A$};
\node (e) at (5.5,1) {$A$};
\node (f) at (5,2.33) {$A$};
\node (g) at (5,3.67) {$W$};
\node (h) at (6.5,3.67) {$A^\vee$};
\node[sdot] (m) at (5,1.67) {};
\node[sdot] (n) at (5,3) {};
\node[right] at (m) {$\;m$};
\node[right] at (n) {$\;n$};
\draw[thick] (a) -- (phi);
\draw[thick] (b) -- (phi);
\draw[thick] (c) -- (phi);
\draw[thick] (d) -- (4.5,1) -- (m) -- (f) -- (n) -- (g);
\draw[thick] (m) -- (e) to[out=-90,in=-90] (6.5,0.8) -- (h);
\end{tikzpicture}
\]
is invertible.
\end{definition}

Note that if $W=\ti$ then we recover the more restrictive definition of a Frobenius structure given in \cite{fs} and \cite{street}.  We call this a \emph{classical Frobenius structure}.

\begin{remark}\label{rmk:rightAhom}
If $A$ is an algebra object then the left dual $A^\vee$ is a categorical right $A$-module with the usual action map
\[ A^\vee\otimes A \arrr{\id_{A^\vee\otimes A}\otimes c_A} A^\vee\otimes A\otimes A\otimes A^\vee \arrr{\id_{A^\vee}\otimes m\otimes\id_{A^\vee}} A^\vee\otimes A\otimes A^\vee\arrr{e_A\otimes\id_{A^\vee}}A^\vee \]
and applying $W\otimes-$ gives the structure of a categorical right $A$-module to $W\otimes A^\vee$.  Then applying associativity of $m$ and using a triangle identity shows that $\varphi$ is automatically a morphism of right $A$-modules.
\end{remark}

The next result shows that our more general definition gives nothing new in the classical situation of finite dimensional algebras over vector spaces.
\begin{proposition}
If $(W,n)$ is a Frobenius structure on an algebra $A$ in $\Vect$ then $W\cong\ti=\kk$.
\end{proposition}
\begin{proof}
As $A$ is dualizable it has finite dimension, and $\dim_\kk A=\dim_\kk A^\vee$.  So, as the map $A\to W\otimes A^\vee$ is invertible, we must have $\dim_\kk W=1$.  So $W\cong\kk$.
\end{proof}

An object $X\in\Ccat$ is called \emph{invertible} if it is left dualizable and $e_X$ and $c_X$ are both isomorphisms.  Then $X^\vee$ is called the \emph{inverse} object of $X$.  Note that the inverse of the triangle identities for $X$ show that $X^\vee$ has left dual $X$ with $e_{X^\vee}=c_X^{-1}$ and $c_{X^\vee}=e_X^{-1}$.  Therefore ${^\vee X}\cong X^\vee$.

We now give some results on Frobenius structures $(W,n)$ under the assumption that $W$ is invertible.

\begin{lemma}\label{lem:adual}
Suppose $A$ has Frobenius structure $(W,n)$ where $W$ has inverse object $V$.  Then $A^\vee\cong V\otimes A$ and composing $\varphi$ with this isomorphism gives $c_W\otimes\id_A:A\to W\otimes V\otimes A$.
\end{lemma}
\begin{proof}
We have an isomorphism $e_W:V\otimes W\arr\sim \ti$, so we construct a map
\[ f: V\otimes A\arrr{\id_V\otimes\varphi} V\otimes W\otimes A^\vee\arrr{e_W\otimes\id_{A^\vee}}A^\vee\]
which, using $e_W^{-1}=c_V$, is an isomorphism with inverse
\[ f^{-1}: A^\vee \arrr{c_V\otimes\id_{A^\vee}} V\otimes W\otimes A^\vee \arrr{\id_V\otimes\varphi^{-1}} V\otimes A.\]
Thus $A^\vee\cong V\otimes A$.  Consider the composition
\[ \varphi_\text{new}:A\arrr\varphi W\otimes A^\vee \arrr{\id_W\otimes f^{-1}} W\otimes V\otimes A. \]
Using naturality of $\otimes$ and the triangle identities, we have:
\begin{align*}
(e_V\otimes\id_A)\circ \varphi_\text{new}
 &= (e_V\otimes\id_A)\circ (\id_W\otimes \id_V\otimes \varphi^{-1})\circ (\id_W\otimes c_V\otimes\id_{A^\vee} )\circ\varphi \\
 &= \varphi^{-1}\circ (e_V\otimes\id_W\otimes \id_{A^\vee})\circ (\id_W\otimes  c_V\otimes\id_{A^\vee} )\circ\varphi \\
 &= \varphi^{-1}\circ (\id_W\otimes \id_{A^\vee})\circ\varphi = \id_A 
\end{align*}
So $\varphi_\text{new}=(e_V\otimes\id_A)^{-1}=c_W\otimes\id_A$.
\end{proof}

\begin{lemma}
Suppose $W$ has inverse object $V$.  
Then the following are equivalent:
\begin{enumerate}[(i)]
\item $A$ is left dualizable and $(A,u,m)$ has Frobenius structure $(W,n)$;
\item $A$ is right dualizable and the map
\[ \varphi':A \larr{c^r_A\otimes\id_A} {^\vee A}\otimes A\otimes A\larr{\id_{^\vee\!A}\otimes\kappa} {^\vee A}\otimes W \]
is invertible.
\end{enumerate}
If the equivalent conditions hold, then $A^\vee\cong V\otimes A$ and ${^\vee A}\cong A\otimes V$.
\end{lemma}
\begin{proof}
The following diagrams illustrate the constructions used below:
\[
\begin{tikzpicture}
[auto, block/.style ={rectangle, draw=black, thick, 
text width=2em, align=center, rounded corners, minimum height=2em},
 wblock/.style ={rectangle, draw=black, thick, 
 text width=2em, align=center, rounded corners, minimum height=2em, minimum width=4em},
 sdot/.style = {circle,draw=black, minimum width=5pt, fill=black, inner sep=0pt]} ]
 
\node[wblock] (t1) at (0.55,3) {};
\node (t1a) at (0.2,3) {$A$};
\node (t1b) at (0.9,3) {$V$};
\node[block] (t2) at (2.5,3) {$A$};
\node[sdot] (phiinv) at (2.5,2.1) {};
\node[right] at (phiinv) {$\;\varphi^{-1}$};

\coordinate (m1) at (0.2,1.3) {};
\coordinate (m2) at (0.9,1.3) {};
\node (m3) at (2,1.5) {$W$};
\node (m4) at (3,1.5) {$A^\vee$};
\node (m5) at (5,1.5) {$W$};
\node (m6) at (6,1.5) {$A^\vee$};
\coordinate (m7) at (7.2,1.7) {};
\coordinate (m8) at (7.9,1.7) {};

\node[wblock] (b2) at (7.55,0) {};
\node (b2a) at (7.2,0) {$A$};
\node (b2b) at (7.9,0) {$V$};
\node[block] (b1) at (5.5,0) {$A$};
\node[sdot] (phi) at (5.5,0.9) {};
\node[right] at (phi) {$\;\varphi$};

\draw[thick] (b1) -- (phi);
\draw[thick] (m5) -- (phi);
\draw[thick] (m6) -- (phi);

\draw[thick] (t2) -- (phiinv);
\draw[thick] (m3) -- (phiinv);
\draw[thick] (m4) -- (phiinv);

\draw[thick] (0.2,2.65) to[out=-90,in=90]  (m1) to[out=-90,in=-90] (m4);
\draw[thick] (7.9,0.35) to[out=90,in=-90]  (m8) to[out=90,in=90] (m5);

\draw[thick] (0.9,2.65) to[out=-90,in=90]  (m2) to[out=-90,in=-90] (m3);
\draw[thick] (7.2,0.35) to[out=90,in=-90]  (m7) to[out=90,in=90] (m6);

\node (x1) at (-2.5,3) {${^\vee A}$};
\node (x2) at (-1.5,3) {$A$};
\draw[thick] (x1) to[out=-90,in=-90]  (x2);
\node at (-0.8,2.5) {$=$};

\node (y2) at (10.5,0) {${^\vee A}$};
\node (y1) at (9.5,0) {$A$};
\draw[thick] (y1) to[out=90,in=90]  (y2);
\node at (8.8,0.5) {$=$};
\end{tikzpicture}
\]

Suppose \emph{(i)} holds.  We know that 
$A^\vee\cong V\otimes A$
by Lemma \ref{lem:adual}.  First we define ${^\vee A}= A\otimes V$ and show that it is a right dual of $A$.  
We construct the evaluation as
\[ e^r_A:A\otimes {^\vee A} =A\otimes A\otimes V\arrr{\varphi\otimes\id_{A\otimes V}} W\otimes A^\vee \otimes A\otimes V
\arrr{\id\otimes e_A\otimes \id}W\otimes  V\arrr{e_V}\ti
\]
and the coevaluation as
\[ c^r_A:\ti\arrr{c_A}A\otimes A^\vee\arrr{\id\otimes c_V\otimes \id}
A\otimes V\otimes W\otimes A^\vee\arrr{\id\otimes\varphi^{-1}}
A\otimes V\otimes A={^\vee A} \otimes A\]
and it is straightforward to check that they satisfy the triangle identities.  So we have constructed ${^\vee A}$ explicitly and thus $A$ is right dualizable.

Now we show that $\varphi'$ is invertible.  The maps $\varphi$ and $\varphi'$ are related by  
\[\varphi'=(\id_{^\vee A\otimes W}\otimes e_A) \circ(\id_{^\vee A}\otimes\varphi\otimes \id_A)\circ(c^r_A\otimes\id_A).\]   
Using the description of $c^r_A$ above we deduce that $\varphi'=\id_A\otimes c_V:A\to A\otimes V\otimes W$, giving a description of $\varphi'$ similar to that of Lemma \ref{lem:adual} for $\varphi$.  Then because $V$ is the inverse of $W$ we know that $c_V$ is invertible, and therefore $\varphi'$ is invertible.

For $\emph{(ii)}\then\emph{(i)}$ we can apply $\emph{(i)}\then\emph{(ii)}$ in the monoidal category with opposite tensor product.
\end{proof}

Sometimes the following less symmetric definition is more convenient to work with.  
\begin{definition}\label{def:left-frob}
Suppose $A$ is left dualizable.  A \emph{left Frobenius structure} on the algebra $(A,u,m)$ is a pair $(V,n^\ell)$ where $V\in\Ccat$ is a right dualizable object and
\[n^\ell:V\otimes A\to \ti\] 
is a map in $\Ccat$ such that
\[ \varphi^\ell:V\otimes A\larr{\id_V\otimes \id_A\otimes c_A} V\otimes A\otimes A\otimes A^\vee \larr{\id_{V}\otimes m\otimes \id_{A^\vee}}V\otimes A\otimes A^\vee\larr{n^\ell\otimes \id_{A^\vee}} A^\vee \]
is invertible in $\Ccat$.
\end{definition}
The following result shows that Definitions \ref{def:frob} and \ref{def:left-frob} are equivalent.
\begin{lemma}\label{lem:can-use-left}
Let $(A,u,m)$ be an algebra in $\Ccat$.  Suppose $W\in\Ccat$ is invertible with inverse $V$ and consider a map $n:A\to W$.  Then $(W,n)$ is a twisted Frobenius structure if and only if $(V,n^\ell)$ is a left Frobenius structure, where $n^\ell= e_W\circ(\id_V\otimes n):V\otimes A\to V\otimes W\to \ti$.
\end{lemma}
\begin{proof}
We have $n=(\id_W\otimes n^\ell)\circ(c_W\otimes\id_A)$, so $\varphi^\ell=(e_W\otimes \id_{A^\vee})\circ(\id_V\otimes \varphi):V\otimes A\to V\otimes W\otimes A^\vee\to A^\vee$ and $\varphi=(\id_W\otimes \varphi^\ell)\circ (c_W\otimes \id_A):A\to W\otimes V\otimes A\to W\otimes A^\vee$.  

Suppose $\varphi$ has inverse $\psi$ and define $\psi^\ell=(\id_V\otimes \psi) 
\circ (e_W^{-1}\otimes\id_{A^\vee}):
A^\vee \to V\otimes W\otimes A^\vee \to V\otimes A^\vee$.  
A straightforward check shows $\psi^\ell$ is a two-sided inverse for $\varphi^\ell$.  Conversely, given $\psi^\ell$ inverse to $\varphi^\ell$, construct an inverse $\psi= (c_W^{-1}\otimes\id_{A})\circ (\id_W\otimes \psi^\ell)$  to $\varphi$.
\end{proof}
One could also define a right Frobenius structure with a map $n^r:A\otimes V\to \ti$, and a result analogous to Lemma \ref{lem:can-use-left} holds.


\subsection{Nakayama automorphisms}

Classically, there are at least two ways to define the Nakayama automorphism.  The first is Nakayama's original definition: given a Frobenius $\kk$-algebra $A$ with Frobenius form $n:A\to \kk$, Nakayama showed that there exists an automorphism $\alpha:A\to A$ such that, for all $x,y\in A$, we have $n(\alpha(x)y)=n(yx)$ \cite[Theorem 1]{na2}, now known as the Nakayama automorphism.  One could mimic this definition for a Frobenius algebra object in a braided monoidal category, but we will have reason to consider tensor categories which do not admit a braiding, as the following example shows.

\begin{example}\label{eg:no-braiding}
Let $S=\kk\times\kk$ and write $e_1=(1,0)$ and $e_2=(0,1)$.  Let $\Ccat$ be the category of $S\da S$-bimodules, with $\otimes=\otimes_S$.  Let $M=\spn_\kk\{m\}$ and $N=\spn_\kk\{n\}$ be $1$-dimensional $\kk$-vector spaces with $S\da S$-bimodule structures determined by $e_1me_1=m$ and $e_1ne_2=n$.  Note that $e_2m=0$.  Then $M\otimes_S N\cong N\neq0$ and $N\otimes_SM=0$, so $\Ccat$ cannot be braided.
\end{example}

Composing the multiplication $m:A\otimes_\kk A\to A$ with the pairing $n:A\to \kk$ gives a bilinear form $\kappa:A\otimes_\kk A\to \kk$.  This induces an isomorphism $\kappa^r:A\arr\sim A^*$ of right $A$-modules sending $y\in A$ to the linear functional $x\mapsto \kappa(y,x)$.  One checks that $\kappa^r$ is in fact an isomorphism of $A\da A$-bimodules if one twists the left action on $A^*$ by the inverse Nakayama automorphism $\alpha^{-1}$, and this gives the second definition of $\alpha$.

We can construct the Nakayama automorphism directly as the composition 
\[A\arrr{\kappa^r} A^* \arrr{((\kappa^r)^*)^{-1}}A^{**}\arrr{\text{ev}^{-1}}A.\]
One could mimic this definition for a Frobenius algebra object in a rigid monoidal category.  We will generalise this approach to our setting in Definition \ref{def:nak} below.

Let $A\in\Ccat$.  Throughout, we assume that $(A,u,m)$ has Frobenius structure $(W,n)$ where $W$ is invertible with $W^\vee=V$.  So we have an isomorphism $\varphi: A\to W\otimes A^\vee$ which induces $A^\vee\cong V\otimes A$.

The category $\Ccat$ is pivotal if there is an isomorphism of monoidal functors $\iota:-\to -^{\vee\vee}$ from the identity to the double dual: 
see \cite[Section 4.7]{egno}.  We say that an object $A\in\Ccat$ is \emph{pivotal} if the full subcategory on the object $A$ is pivotal.

The isomorphism $A^\vee\cong V\otimes A$ induces isomorphisms
\[ A^{\vee\vee}\cong (V\otimes A)^\vee \cong A^\vee\otimes V^\vee\cong V\otimes A \otimes V^\vee \]
so if $A$ is pivotal we should get an isomorphism between $V\otimes A$ and $A\otimes V$.  Even though this is not an automorphism, it will play the role of the Nakayama automorphism in our setup.
\begin{definition}\label{def:nak}
Suppose $A$ is pivotal.
The \emph{Nakayama morphism} of $A$ is the map
\[ \alpha: V\otimes A \to A\otimes V \]
defined by
\[ V\otimes A \arrr{\id_V\otimes\varphi} V\otimes W\otimes A^\vee \arrr{e_W\otimes (\varphi^\vee)^{-1}} A^{\vee\vee}\otimes V \arr{\iota_A^{-1}\otimes\id_V} A\otimes V.\]
\end{definition}
\[
\begin{tikzpicture}
[auto, block/.style ={rectangle, draw=black, thick, fill=blue!10, text width=2em, align=center, rounded corners, minimum height=2em},
 sdot/.style = {circle,draw=black, minimum width=5pt, fill=black, inner sep=0pt]} ]
\node (a) at (-2,1.2) {$V$};
\node (b) at (-1,1.2) {$A$};
\node (c) at (-2,2.4) {$A$};
\node (d) at (-1,2.4) {$V$};
\node[sdot] (alpha) at (-1.5,1.8) {};
\node[right] at (alpha) {$\;\alpha$};

\draw[thick] (a) -- (alpha) -- (c);
\draw[thick] (b) -- (alpha) -- (d);

\node at (0,1.8) {$=$};

\node (e) at (1,0) {$V$};
\node (f) at (2.5,0) {$A$};
\node[sdot] (phi) at (2.5,0.6) {};
\node[right] at (phi) {$\;\varphi$};

\node (g) at (2,1.2) {$W$};
\coordinate (i) at (1,1.45) {};
\node (h) at (3,1.2) {$A^\vee$};
\node[sdot] (j) at (3,1.8) {};
\node[right] at (j) {$\;(\varphi^\vee)^{-1}$};

\node (k) at (2.5,2.4) {$A^{\vee\vee}$};
\coordinate (l) at (3.5,2.2) {};
\node[sdot] (iot) at (2.5,3) {};
\node[right] at (iot) {$\;\iota_A^{-1}$};
\node (m) at (2.5,3.6) {$A$};
\node (n) at (3.5,3.6) {$V$};

\draw[thick] (e) to[out=90,in=-90] (i) to[out=90,in=90] (g);
\draw[thick] (g) -- (phi);
\draw[thick] (f) -- (phi) -- (h) -- (j) -- (k) -- (iot) -- (m);
\draw[thick] (n) to[out=-90,in=90] (l) -- (j);

\end{tikzpicture}
\]

By abuse of notation, we write $\alpha^2$ for the map
\[ \alpha^2:V\otimes V\otimes A\arr{\id_V\otimes\alpha} V\otimes A\otimes V\arr{\alpha\otimes\id_V}  A\otimes V\otimes V.\]
Similarly, for any $n\geq1$, let $\alpha^n$ denote the map
\[ \alpha^n: V^{\otimes n}\otimes A\to A\otimes V^{\otimes n}.\]
\begin{definition}\label{def:nak-order}
We write $\alpha^n=1$ if there exists an invertible map $\xi:V^{\otimes n}\arr\sim \ti$ in $\Ccat$ such that the composition
\[ A\arr\sim \ti\otimes A\arr{\xi^{-1}\otimes\id_A} 
V^{\otimes n}\otimes A \arr{\alpha^n}  A\otimes V^{\otimes n}
\arr{\id_A\otimes\xi}  A\otimes \ti\arr\sim A \]
is the identity on $A$.  We say $\alpha$ has order $n\in\{1,2,\ldots,\infty\}$ if $n$ is minimal such that $\alpha^n=1$. 
\end{definition}


\subsection{Transfer of structure}

Let $\Ccat$, $\D$ be monoidal categories.
Recall that a monoidal functor from $\Ccat$ to $\D$ is a pair $(F,J)$ consisting of a functor $F:\Ccat\to\D$ and isomorphisms $J_{C,C'}:F(C)\otimes F(C')\to F(C\otimes C')$ natural in both $C,C'\in\Ccat$.  It should satisfy unitality and associativity conditions.  Note that monoidal functors respect duals \cite[Exercise 2.10.6]{egno}.

Fix a monoidal functor $(F,J):\Ccat\to\D$ and an algebra $(A,m,u)$ in $\Ccat$.  Then we get an algebra $(B=F(A),m_B,u_B)$ in $\D$ with structure maps defined using $F$ and $J$.  
\begin{proposition}\label{prop:monfun}
Any monoidal functor $(F,J):\Ccat\to\D$ sends a Frobenius structure $(W,n)$ on $(A,m,u)$ to a Frobenius structure $(Y,n_B)=(F(W),F(n))$ on $(B,m_B,u_B)$.
\end{proposition}
\begin{proof}
We need to show that the map $\varphi_B:B\to Y\otimes B^\vee$ is invertible.  Applying the definitions, and using the unitality, naturality, and associativity conditions on $J$, we get that the following diagram commutes:
\[\xymatrix @C=40pt { F(A) \ar[r]^-{F(\varphi_A)} & F(W\otimes A^\vee) \\
B \ar@{=}[u] \ar[r]_-{\varphi_B} & Y\otimes B^\vee \ar[u]_{J_{W,A^\vee}} }\]
As $F(\varphi_A)$ and $J$ are invertible, so is $\varphi_B$.
\end{proof}

Keep the notation from above, and let $V\in\Ccat$ and $X\in\D$ be the left duals of $W$ and $Y$, respectively.
\begin{lemma}\label{lem:nak-transf}
Suppose $A\in \Ccat$ and $B\in \D$ are both pivotal, with $\iota_B^\D=F(\iota_A^\Ccat)$.  Let $\alpha_A$ and $\alpha_B$ denote the Nakayama morphisms of $A$ and $B$.  Then the following diagram commutes:
\[\xymatrix @C=40pt { F(V\otimes A) \ar[r]^-{F(\varphi_A)} & F(A\otimes V) \\
X\otimes B \ar[u]^{J_{V,A}} \ar[r]_-{\varphi_B} & B\otimes X \ar[u]_{J_{A,V}} }\]
\end{lemma}
\begin{proof}
This is another straightforward check, using conditions on $J$ and compatibility of the pivotal structures.
\end{proof}
Keep the conditions and notation of Lemma \ref{lem:nak-transf}.
Applying the same methods as above, together with naturality of tensor products, gives the following result:
\begin{corollary}
Suppose $\alpha_A^n=1$ and $\iota_B^\D=F(\iota_A^\Ccat)$.  Then $\alpha_B^n=1$.
\end{corollary}

\subsection{Examples}
Fix a field $\kk$.
\begin{example}\label{ex:toyTL}
Let $\Ccat$ be a semisimple $\kk$-category with two non-isomorphic objects, $\ti$ and $X$, and all their finite direct sums.  Define a strict monoidal structure on $\Ccat$ by $X\otimes X=\ti$.  This structure is rigid, with $X^\vee=X={^\vee X}$.  Let $A=\ti\oplus X$.  Define $u:\ti\to A$ by inclusion and $m:A\otimes A\to A$ by inclusion for $\ti\otimes\ti$, $\ti\otimes X$, and $X\otimes \ti$, and by the zero map for $X\otimes X$.  Then $(A,m,u)$ is an algebra in $\Ccat$.

Define a map $n:A\to X$ by projection, and note that $\kappa=n\circ m:A\otimes A\to X$ acts as the identity on $\ti\otimes X$ and $X\otimes \ti$, and as zero on $\ti\otimes\ti$ and $X\otimes X$.  We claim that $(X,n)$ is a Frobenius structure.  By properties of duals of direct sums, the map $c_A:\ti\to A\otimes A$ has identity components to $\ti\otimes\ti$ and $X\otimes X$ and zero components to $\ti\otimes X$ and $X\otimes \ti$.  Therefore the map $\varphi:A\arrr{\id\otimes c}A\otimes A\otimes A\arrr{\kappa\otimes\ti} X\otimes A$ is defined on each summand by the following diagram, where all arrows are isomorphisms:
\[\xymatrix @R=0pt{
  & \ti\otimes \ti\otimes \ti &  \\
\ti \ar[ur]!<-3ex,1ex>\ar[ddr]!<-3ex,-1ex>  & \ti\otimes \ti\otimes X & X\otimes \ti \cong X \\
  & \ti\otimes X\otimes \ti \ar!<0ex,-1.5ex>;[ur]!<0ex,2ex>  &  \\
  & \ti\otimes X\otimes X \ar!<1ex,3.5ex>;[ddr]!<-4ex,0ex> &  \\
  & X\otimes \ti\otimes \ti \ar!<3ex,-2.5ex>;[uuur]!<-4ex,0ex> &  \\
X \ar[ur]!<-3ex,1ex>\ar[ddr]!<-3ex,-1ex> & X\otimes \ti\otimes X \ar[r] & X\otimes X \cong \ti \\
  & X\otimes X\otimes \ti &  \\
  & X\otimes X\otimes X &  
}\]
We see that $\varphi$ is an isomorphism having summands $\ti\arr\sim X\otimes X$ and $X\arr\sim X\otimes \ti$.  Therefore the dual map $\varphi^\vee$ has summands $X\otimes X\arr\sim \ti$ and $\ti\otimes X\arr\sim X$.

Consider the trivial pivotal structure, where we identify an object with its double dual using the identity map.  Then the Nakayama morphism $X\otimes A\to A\otimes X$ is the direct sum of the following isomorphisms:
\[\xymatrix @R=0pt{
  X\otimes \ti \ar[r] & X\otimes X\otimes X \ar[r] & \ti\otimes X \ar[r]& \ti\otimes X \\
  X\otimes X \ar[r] & X\otimes X\otimes \ti \ar[r] & X\otimes X \ar[r]& X\otimes X
}\]
\end{example}

\begin{example}\label{ex:toyPi}
Let $Q$ be the following quiver
\[
\xymatrix{
1 \ar@/^/[r]^{a} & 2\ar@/^/[l]^{b}
}
\]
and let $A=\kk Q/(ab,ba)$ be its path algebra modulo paths of length $2$, so $A$ has basis $\{e_1,e_2,a,b\}$.  Then $A$ is a Frobenius $\kk$-algebra with Frobenius form $n:A\to \kk$ defined by $n(a)=n(b)=1$ and $n(e_1)=n(e_2)=0$.  Its (classical) Nakayama automorphism is the map of $\kk$-vector spaces $A\to A$ which interchanges $e_1\leftrightarrow e_2$ and $a\leftrightarrow b$.

Let $S$ denote the subalgebra generated by $e_1$ and $e_2$ and let $\D=S\mModm S$ be the category of $S\da S$-bimodules, as in Example \ref{eg:no-braiding}.  Note that $A$ has left and right $S$-actions and the multiplication is balanced over $S$, so $A$ is an algebra object in $\D$.  However, all maps $A\to S$ in $\D$ must send $a$ and $b$ to zero, so one can deduce that there are no classical Frobenius structures on the algebra $A\in\D$.

Let $W$ be a $2$-dimensional vector space with basis $\{w_{12}, w_{21}\}$.  We give $W$ an $S\da S$-bimodule structure by $e_iw_{ij}e_j=w_{ij}$.  Note that $W$ is self-inverse: $W\otimes W\cong S$, and therefore $W^\vee=W$. Define $n:A\to W$ by $n(a)=w_{12}$, $n(b)=w_{21}$, and $n(e_1)=n(e_2)=0$.  Then one can check that $(W,n)$ is a Frobenius structure on $A\in\D$.  Its Nakayama morphism is defined by:
\begin{align*}
\alpha:W\otimes A &\to A\otimes W\\
w_{12}\otimes e_2 &\mapsto  e_1\otimes w_{12}\\
w_{21}\otimes e_1 &\mapsto  e_2\otimes w_{21}\\
w_{12}\otimes b &\mapsto  a\otimes w_{21}\\
w_{21}\otimes a &\mapsto  b\otimes w_{12}\\
\end{align*}
On ignoring the $w_{ij}$ factors of the tensor product, one can see the same behaviour that appears in the (classical) Nakayama automorphism of the $\kk$-algebra.

One can check that $A$ is a $\beta$-Frobenius extension of $S$ (as defined below), where $\beta$ is the automorphism of $S$ which interchanges $e_1\leftrightarrow e_2$, and so this example illustrates Proposition \ref{prop:2ndkind} below.  Also, there is a monoidal functor $\Ccat\to\D$ from the monoidal category in Example \ref{ex:toyTL} to the monoidal category here, sending $\ti$ to $S$ and $X$ to the $S\da S$-bimodule with $\kk$-vector space basis $\{a,b\}$, so this example illustrates Proposition \ref{prop:monfun} above.
\end{example}

\begin{remark}
Our main example, involving the Temperley-Lieb category, will be given in Section \ref{s:TL}.
When $\kk=\Cnum$, Examples \ref{ex:toyTL} and \ref{ex:toyPi} are both baby examples of this main result.
\end{remark}

We now give a class of examples which shows how our construction generalises Frobenius extensions of the second kind.

\begin{definition}[{Nakayama-Tsuzuku \cite[Section 1]{nt}}]
Let $R$ be a ring with subring $S$, and let $\beta:S\to S$ be a ring automorphism.  
Let ${S_\beta}$ denote the identity $S\da S$-bimodule with right action twisted by $\beta$.
We say that $R$ is a Frobenius extension of the second kind, or a $\beta$-Frobenius extension, of $S$ if:
\begin{enumerate}
\item $R$ is finitely generated and projective as a right $S$-module, and
\item there is an isomorphism $f:R\to\Hom_{\Modm S}(R,{S_\beta})$ of $S\da R$-bimodules.
\end{enumerate}
\end{definition}

Let $\Ccat=S\mModm S$ be the monoidal category of $S\da S$-bimodules, with tensor product over $S$.
Note that $R$ is an algebra object in $\Ccat$.   
\begin{proposition}\label{prop:2ndkind}
$R$ is a $\beta$-Frobenius extension of $S$ if and only if there exists a map $n:R\to S_\beta$ such that $(S_\beta,n)$ is a Frobenius structure on $R$.
\end{proposition}
\begin{proof}
By a standard argument, $M\in S\mModm S$ is left dualizable if and only if it is finitely generated and projective as a right $S$-module.  In this case, the left dual is $M^\vee=\Hom_{\Modm S}(M,S)$, homomorphisms of right $S$-modules.  A statement can be found in \cite[Exercise 2.10.16]{egno} (note the left/right correction in the online errata) and an explanation in \cite{yuan}, where opposite conventions for left and right duals are used.

First suppose that that $(S_\beta,n)$ is a Frobenius structure on $R$.  As $R$ is left dualizable, we know that $R$ is finitely generated and projective as a right $S$-module.  By assumption the map $ \varphi:A \to W\otimes A^\vee$ from Definition \ref{def:frob} is an isomorphism $R\cong S_\beta\otimes_B\Hom_{\Modm S}(R,S)$ of $S\da S$-bimodules, and by Remark \ref{rmk:rightAhom} it is in fact an isomorphism of $S\da R$-bimodules.  Finally, as $R$ is projective as a right $S$-module the natural map $S_\beta\otimes_S\Hom_{\Modm S}(R,S)\to \Hom_{\Modm S}(R,S_\beta)$ of $S\da R$-bimodules is an isomorphism.  

Conversely, suppose $R$ is a $\beta$-Frobenius extension of $S$.  The first condition ensures that $R\in\Ccat$ is left dualizable.  The second condition gives a map $n=f(1_R):R\to S_\beta$.
The proof of \cite[Proposition 4]{nt} shows that $f$ can be reconstructed from $n$: we send $r\in R$ to the homomorphism sending $x\in R$ to $n(rx)$.  Using the identification $S_\beta\otimes_S\Hom_{\Modm S}(R,S)\cong \Hom_{\Modm S}(R,S_\beta)$ this is the formula for $ \varphi:A \to W\otimes A^\vee$, so $\varphi$ is invertible.
\end{proof}


\section{The Temperley-Lieb category}\label{s:TL}

\subsection{Jones-Wenzl projections}\label{ss:jw}
Let $\N$ denote the natural numbers, which for us include zero.  We consider planar diagrams $m\to n$ on a rectangle, with $m$ dots on the bottom and $n$ dots on the top, and non-crossing lines so that every dot is joined to exactly one line.  We allow finitely many loops, which are not joined to any dots.

Recall that the idempotent completion of a category has pairs $(f,x)$ as objects, where $x$ is an object in our original category and $f:x\to x$ is an idempotent.  Maps $(f,x)\to (g,y)$ are of the form $ g\varphi f$, where $\varphi:x\to y$ is a map in our original category, and we sometimes just write them as $\phi:(f,x)\to (g,y)$.

We define the following structures:
\begin{itemize}
\item $\D$ is the category with $\ob \D=\N$ and 
\[ \D(m,n)=\begin{cases}
\emptyset &\text{if } 2 \nmid m+n;\\
\text{diagrams $m\to n$ up to isotopy} &\text{otherwise}.\end{cases}\]
Composition is by vertical stacking, and we have a monoidal structure where $m_1\otimes m_2=m_1+m_2$ and morphisms are tensored by joining diagrams horizontally.  Note that $\D$ is the free pivotal monoidal category on one self-dual generator: see \cite[Proposition 1.1]{ab}.

\item Let $\Cnum \D$ denote the linearization of $\D$, so $\ob \Cnum \D=\ob \D$ and $\Cnum \D(m,n)$ is the $\Cnum$-vector space with basis $\D(m,n)$.  Composition in $\Cnum \D$ is defined by linearity.  
\item Fix $k\geq1$.  Let $t=e^{\frac{\pi i}{2k+4}}$ be a fixed $4k+8$th root of unity and let $q=t^2=e^{\frac{\pi i}{k+2}}$.  
\item Let $[i]=\frac{q^i-q^{-i}}{q-q^{-1}}\in\R$ denote the $i$th quantum integer and write $\delta=-[2]=-q-q^{-1}$.  
\item Define $TL_k$ to be the additive and idempotent completion of the quotient of $\Cnum T$ by the relation that every loop is equal to the constant $\delta$.  We sometimes write this as $TL$, dropping the subscript $k$.  The category $TL$ inherits a monoidal structure from $\D$.
\end{itemize}

Note that $q^{k+2}=-1$, so $[k+2-j]=[j]$.  In particular, $[k+2]=[0]=0$.

Let $\id_1:1\to 1$, $\cup:0\to 2$, and $\cap:2\to 0$ denote the unique such morphisms in $\D$, then define 
\[  \cup_{i,n}=\id_1^{\otimes i} \otimes \cup \otimes \id_1^{\otimes n-i}:n\to n+2 \]
and
\[  \cap_{i,n}=\id_1^{\otimes i} \otimes \cap \otimes \id_1^{\otimes n-i}:n+2\to n \]
both for $0\leq i\leq n$.  
Let $U_{i,n}=\cup_{i,n-2}\circ\cap_{i,n-2}:n\to n$ denote the cap-cup morphism for $1\leq i\leq n-2$, so $U_i^2=\delta U_i$ in $TL$.  Define iterated cups $\Cup_{i,m,n}:n\to n+2m$ and caps $\Cap_{i,m,n}:n+2m\to n$ by
\[  \Cup_{i,1,n}= \cup_{i,n} \;\; \text{ and } \;\;  \Cup_{i,m,n}=\Cup_{i+1,m-1,n+2}\circ \cup_{i,n} \]
and
\[  \Cap_{i,1,n}= \cup_{i,n} \;\; \text{ and } \;\;  \Cap_{i,m,n}=\cap_{i,n}\circ  \Cap_{i+1,m-1,n+2}\]
both for $0\leq i\leq n$ and $m\geq 1$.  These morphisms in $\D$ all determine morphisms in $TL$, which we denote by the same symbols, via the canonical functors $\D\into \Cnum\D \to TL$.

Following \cite{wenzl87}, after \cite{jones}, we have Jones-Wenzl projections $f_i:i\to i$ in $TL$, for $0\leq i\leq k+1$, defined recursively by $f_0=\id_0$, 
$f_1=\id_1$, and 
\[f_i=f_{i-1}\otimes 1+\frac{[i-1]}{[i]} (f_{i-1}\otimes 1)\circ U_{i-1}\circ (f_{i-1}\otimes 1).\]
We record some useful facts about the Jones-Wenzl projections.  Pictures follow after the proof.
\begin{proposition}\label{prop:jw-facts}
We have:
\begin{enumerate}[(i)]
\item $f_i^2=f_i$,
\item $f_n\circ \cup_{i,n-2}=0$ for any $0\leq i\leq n-2$,
\item $\cap_{i,n-2}\circ f_n=0$ for any $0\leq i\leq n-2$,
\item $\cap_{n,n}\circ (f_{n+1}\otimes\id_1)\circ \cup_{n,n}=\dfrac{-[n+2]}{[n+1]}f_n$,
\item $\Cap_{n,j,n}\circ (f_{n+j}\otimes\id_1^{\otimes j})\circ \Cup_{n,j,n}=(-1)^j\dfrac{[n+j+1]}{[n+1]}f_n$.
\end{enumerate}
\end{proposition}
\begin{proof}
For the first four properties, see \cite[Chapter XII]{turaev}, especially Exercise 4.6.  A more detailed explanation, with different sign conventions, can be found in \cite[Theorem 2.3.2]{chen}.  The fifth property follows from the fourth by induction.  
\end{proof}
We illustrate the previous proposition graphically, with maps drawn as green boxes.  Our pictures should be read going upwards, and we sometimes represent $j$ strands as a single strand from the object $j$.
\[
\begin{tikzpicture}
[auto,
 block/.style ={rectangle, draw=black, thick, fill=green!10, text width=2em, align=center, rounded corners, minimum height=2em}
]
\node (n1) at (-1,1.5) {\textit{(i)}};
\node[block] (p1) at (0,2.1) {$f_i$};
\node[block] (p2) at (0,0.9) {$f_i$};
\node[block] (p3) at (2,1.5) {$f_i$};
\node (t1) at (0,3) {$i$};
\node (t2) at (2,3) {$i$};
\node (b1) at (0,0) {$i$};
\node (b2) at (2,0) {$i$};
\node (s1) at (1,1.5) {$=$};
\draw[thick] (t1) -- (p1) -- (p2) -- (b1);
\draw[thick] (t2) -- (p3) -- (b2);
\node[block] (p1) at (5,1.5) {$f_n$};
\node (t1) at (5,2.5) {$n$};
\node (b1) at (4.2,0.5) {$i$};
\node (b2) at (5.8,0.5) {$n-i-2$};
\node (s1) at (6,1.5) {$=$};
\node (c1) at (6.6,1.52) {$0$};
\coordinate (i1) at (5,0.8) {};
\node (n2) at (3.8,1.5) {\textit{(ii)}};
\draw[thick] (t1) -- (p1);
\draw[thick] (p1) -- (b1);
\draw[thick] (p1) -- (b2);
\draw[thick] (p1) to[out=-110, in=180] (i1) to[out=0, in=-70]  (p1);
\node (n3) at (8,1.5) {\textit{(iii)}};
\node[block] (p1) at (9.5,1.5) {$f_n$};
\node (t1) at (9.5,0.5) {$n$};
\node (b1) at (8.7,2.5) {$i$};
\node (b2) at (10.3,2.5) {$n-i-2$};
\node (s1) at (10.5,1.5) {$=$};
\node (c1) at (11.1,1.5) {$0$};
\coordinate (i1) at (9.5,2.2) {};
\draw[thick] (t1) -- (p1);
\draw[thick] (p1) -- (b1);
\draw[thick] (p1) -- (b2);
\draw[thick] (p1) to[out=110, in=180] (i1) to[out=0, in=70]  (p1);
\end{tikzpicture}
\]
\[
\begin{tikzpicture}
[auto,
 block/.style ={rectangle, draw=black, thick, fill=green!10, text width=2em, align=center, rounded corners, minimum height=2em}
]
\node (n4) at (-0.5,1) {\textit{(iv)}};
\node[block] (p1) at (0.5,1) {$f_{n+1}$};
\node[block] (p2) at (4,1) {$f_{n}$};
\node (t1) at (0.5,2) {$n$};
\node (t2) at (4,2) {$n$};
\node (b1) at (0.5,0) {$n$};
\node (b2) at (4,0) {$n$};
\node (s1) at (1.75,1) {$=$};
\node (c1) at (2.7,1) {$\dfrac{-[n+2]}{[n+1]}$};
\coordinate (i1) at (1.3,1) {};
\draw[thick, transform canvas={xshift=-1ex}] (t1) -- (p1) -- (b1);
\draw[thick] (t2) -- (p2) -- (b2);
\draw[thick] (p1) to[out=60,in=90] (i1) to[out=-90,in=-60] (p1);
\node (n5) at (6,1) {\textit{(v)}};
\node[block] (p1) at (7,1) {$f_{n+j}$};
\node (p2) at (8.2,1) {$j$};
\node[block] (p3) at (12.5,1) {$f_{n}$};
\node (t1) at (7,2) {$n$};
\node (t2) at (12.5,2) {$n$};
\node (b1) at (7,0) {$n$};
\node (b2) at (12.5,0) {$n$};
\node (s1) at (8.8,1) {$=$};
\node (c1) at (10.6,1) {$(-1)^j\dfrac{[n+j+1]}{[n+1]}$};
\node (i1) at (7.8,0.9) {};
\draw[thick, transform canvas={xshift=-1ex}] (t1) -- (p1) -- (b1);
\draw[thick] (t2) -- (p3) -- (b2);
\draw[thick] (p1) to[out=60,in=90] (8,1);
\draw[thick] (8,1) to[out=-90,in=-60] (p1);
\end{tikzpicture}
\]

\begin{proposition}[Absorption rules]\label{prop:absorb}
We have:
\begin{enumerate}[(i)]
\item $f_n\circ (\id_{i}\otimes f_j\otimes \id_{n-i-j})=f_n$,
\item $(\id_{i}\otimes f_j\otimes \id_{n-i-j})\circ f_n =f_n$.
\end{enumerate}
\end{proposition}
\begin{proof}
See \cite[Exercise XII.4.6]{turaev} or the proof of part 2 of \cite[Theorem 2.3.2]{chen}.
\end{proof}
Again, we can illustrate these graphically:
\[
\begin{tikzpicture}
[auto,
 block/.style ={rectangle, draw=black, thick, fill=green!10, text width=2em, align=center, rounded corners, minimum height=2em}
]
\node (n1) at (-0.5,0.5) {\textit{(i)}};
\node[block] (p1) at (1,1.1) {$f_n$};
\node[block] (p2) at (1,-0.1) {$f_j$};
\node[block] (p3) at (3.8,0.5) {$f_n$};
\node (t1) at (1,2) {$n$};
\node (t2) at (3.8,2) {$n$};
\node (b1) at (-0.2,-1) {$i$};
\node (b2) at (1,-1) {$j$};
\node (b3) at (2.2,-1) {$n-i-j$};
\node (b4) at (3.8,-1) {$n$};
\node (s1) at (2.5,0.5) {$=$};
\draw[thick] (t1) -- (p1) -- (p2) -- (b2);
\draw[thick] (p1) to[out=-130,in=90] (b1);
\draw[thick] (p1) to[out=-50,in=90] (b3);
\draw[thick] (t2) -- (p3) -- (b4);
\node (n2) at (5.5,0.5) {\textit{(ii)}};
\node[block] (p1) at (7,-0.1) {$f_n$};
\node[block] (p2) at (7,1.1) {$f_j$};
\node[block] (p3) at (9.8,0.5) {$f_n$};
\node (t1) at (7,-1) {$n$};
\node (t2) at (9.8,2) {$n$};
\node (b1) at (5.8,2) {$i$};
\node (b2) at (7,2) {$j$};
\node (b3) at (8.2,2) {$n-i-j$};
\node (b4) at (9.8,-1) {$n$};
\node (s1) at (8.5,0.5) {$=$};
\draw[thick] (t1) -- (p1) -- (p2) -- (b2);
\draw[thick] (p1) to[out=130,in=-90] (b1);
\draw[thick] (p1) to[out=50,in=-90] (b3);
\draw[thick] (t2) -- (p3) -- (b4);
\end{tikzpicture}
\]

As the Jones-Wenzl projections are idempotents, they induce objects $F_i=(i,f_i)$ in $TL$.  Note that $F_0=\ti$ is the unit object.

Let $(F_{k+1})$ denote the ideal of morphisms generated by $F_{k+1}$.
\begin{definition}
The \emph{reduced Temperley-Lieb category}, or \emph{Temperley-Lieb-Jones category} at level $k$, is $\tl=\tlk=TL/(F_{k+1})$.
\end{definition}
The following result is standard: see \cite[page 47 and Theorem 5.4.5]{chen}.
\begin{theorem}\label{thm:tlss}
$\tl$ is semisimple with non-isomorphic simple objects $F_0$, $F_1$, $\ldots$, $F_k$, and
\[ F_1\otimes F_i\cong F_{i-1}\oplus F_{i+1} \]
where $F_\ell=0$ for $\ell>k$.
\end{theorem}
\begin{remark}
In fact, $F_i\otimes F_j$ always contains $F_{i+j}$ as a summand with multiplicity $1$: see \cite[page 47]{chen}.
\end{remark}

Recall the definition of maps in the idempotent completion from the start of Section \ref{ss:jw}.
\begin{definition}\label{def:pi-iota}
Let $i,j\geq0$ with $i+j\leq n$.
We define maps
\[ \iota_{i,j}:F_{i+j}\to F_i\otimes F_j \;\; \text{ and } \;\; \pi_{i,j}:F_i\otimes F_j \to F_{i+j}\]
using $f_{i+j}$ as follows: $\iota_{i,j}=(f_i\otimes f_j)f_{i+j}f_{i+j}$ and $\pi_{i,j}=f_{i+j}f_{i+j}(f_i\otimes f_j)$.
\end{definition}
We draw these morphisms as follows, with Jones-Wenzl projections drawn as blue boxes:
\[
\begin{tikzpicture}
[auto, block/.style ={rectangle, draw=black, thick, fill=blue!10, text width=2em, align=center, rounded corners, minimum height=2em},
 dot/.style = {circle,draw=black, thick, fill=green!10} ]
\node[block] (p1) at (1,0) {$F_{i+j}$};
\node[block] (p2) at (0,2) {$F_i$};
\node[block] (p3) at (2,2) {$F_j$};
\node[dot] (ma) at (1,1) {$\iota$};
\draw[thick] (p1) -- (ma) -- (p2);
\draw[thick] (p3) -- (ma);
\node[block] (p1) at (5,2) {$F_{i+j}$};
\node[block] (p2) at (4,0) {$F_i$};
\node[block] (p3) at (6,0) {$F_j$};
\node[dot] (mb) at (5,1) {$\pi$};
\draw[thick] (p1) -- (mb) -- (p2);
\draw[thick] (p3) -- (mb);
\end{tikzpicture}
\]
Note that $\pi_{i,j}\circ\iota_{i,j}=\id_{F_{i+j}}$ by the absorption rules.  But $\iota_{i,j} \circ\pi_{i,j}\neq\id_{F_i\otimes F_j}$: for example, when $i=j=1$ we have $F_1\otimes F_1\cong F_2\oplus F_0$ and $\iota_{1,1} \circ\pi_{1,1}$ projects to the $F_2$ summand.
\begin{lemma}\label{lem:pi-iota-absorb}
For maps $g:i\to i$ and $h:j\to j$ we have the following equalities of maps in $TL$:
\[ (g\otimes h)\circ \iota_{i,j}=g\otimes h:F_{i+j}\to F_i\otimes F_j \;\; \text{ and }\;\; \pi_{i,j}\circ (g\otimes h)=g\otimes h:F_i\otimes F_j\to F_{i+j}.\]
\end{lemma}
\begin{proof}
From Definition \ref{def:pi-iota} we have
\[(g\otimes h)\circ \iota_{i,j}
=(f_i\otimes f_j)(g\otimes h)(f_i\otimes f_j)f_{i+j}
=(f_i\otimes f_j)(g\otimes h)f_{i+j}
=g\otimes h:F_{i+j}\to F_i\otimes F_j
\]
by Proposition \ref{prop:absorb}; the equation for $\pi$ is similar.
\end{proof}
By Lemma \ref{lem:pi-iota-absorb} with $g$ and $h$ being identity maps, we can represent $\iota$ and $\pi$ as follows:
\[
\begin{tikzpicture}
[auto, block/.style ={rectangle, draw=black, thick, fill=blue!10, text width=2em, align=center, rounded corners, minimum height=2em},
 wblock/.style ={rectangle, draw=black, thick, fill=blue!10, text width=2em, align=center, rounded corners, minimum height=2em, minimum width=4em},
 dot/.style = {circle,draw=black, thick, fill=green!10} ]
\node[wblock] (p1) at (1,0) {$F_{i+j}$};
\node[block] (p2) at (0.4,1.4) {$F_i$};
\node[block] (p3) at (1.6,1.4) {$F_j$};
\node (ca) at (1,1.18) {};
\node (sa) at (-0.5,0.7) {$\pi =$};
\draw[thick, transform canvas={xshift=-2ex}] (ca) -- (p1);
\draw[thick, transform canvas={xshift=2ex}] (ca) -- (p1);
\node[wblock] (p4) at (5,1.4) {$F_{i+j}$};
\node[block] (p5) at (4.4,0) {$F_i$};
\node[block] (p6) at (5.6,0) {$F_j$};
\node (cb) at (5,0.22) {};
\node (sb) at (3.5,0.7) {$\iota =$};
\draw[thick, transform canvas={xshift=-2ex}] (cb) -- (p4);
\draw[thick, transform canvas={xshift=2ex}] (cb) -- (p4);
\end{tikzpicture}
\]

\begin{definition}\label{def:sigma}
Let $\Sigma=\bigoplus_{i=0}^kF_i$.  Let $u:\ti=F_0\into \Sigma$ be the inclusion map.  Let $m=\Sigma\otimes\Sigma\to \Sigma$ be defined by summing maps $m_{i,j}:F_i\otimes F_j\to F_{i+j}\into \Sigma$, where $m_{i,j}=\pi_{i,j}:F_i\otimes F_j\to F_{i+j}$ if $i+j\leq k$, and is zero otherwise.
\end{definition}
The following result is well-known: see, for example, \cite[Section 5.4.3]{cooper}.  For completeness, we include a proof in our context.
\begin{lemma} 
$(\Sigma,u,m)$ is an algebra object in $\tl$.
\end{lemma}
\begin{proof}
We can check everything on graded components.
First we check left and right units.  As the empty diagram $\ti$ is a strict identity, the left and right unit isomorphisms are just identity maps  But the multiplication maps $f_j:F_0F_j\to F_j$ and $f_i:F_iF_0\to F_i$ are identity maps in our idempotent completion.  To check associativity we need commutativity of the following diagram
\[\xymatrix @C=34pt {
F_h\otimes F_i\otimes F_j \ar[r]^{\id\otimes f_{i+j}}\ar[d]^{f_{h+i}\otimes \id} & F_h\otimes F_{i+j} \ar[d]^{f_{h+i+j}}  \\
F_{h+i}\otimes F_j \ar[r]^{f_{h+i+j}} & F_{h+i+j} 
}\]
and this follows from the absorption rules (Proposition \ref{prop:absorb}).
\end{proof}

\begin{lemma}\label{lem:no-classical}
$\Sigma$ does not admit a classical Frobenius algebra structure.
\end{lemma}
\begin{proof}
If $W=\ti$ then, by Theorem \ref{thm:tlss}, the map $n:\Sigma\to \ti$ from Definition \ref{def:frob} can only be nonzero on the summand $F_0$ of $\Sigma=\bigoplus_{i=0}^kF_i$.  Therefore $\varphi$ must act as zero on $\bigoplus_{i>0}^kF_i$ and so cannot be an isomorphism.
\end{proof}

\subsection{Pivotal and braided structures}

\begin{lemma}\label{lem:pull-jw-round}
$\Cap_{0,i,0}(f_i\otimes \id_i)=\Cap_{0,i,0}(\id_i\otimes f_i)$ and
$(f_i\otimes \id_i)\Cup_{0,i,0} = (\id_i\otimes f_i)\Cup_{0,i,0}$.
\end{lemma}
\begin{proof}
This follows from left-right symmetry of the Jones-Wenzl projections: see \cite[Remark 2.3.5]{chen}.
\end{proof}
As a consequence of Lemma \ref{lem:pull-jw-round}, we get:
\begin{lemma}\label{lem:duals}
Every simple object is self-dual, with evaluation $e_i:F_i\otimes F_i\to\ti$ and coevaluation $c_i:\ti\to F_i\otimes F_i$ given by cupping and capping: $e_i=\Cap_{0,i,0}(f_i\otimes f_i)$ and $c_i=(f_i\otimes f_i)\Cup_{0,i,0}$.
\end{lemma}
In particular, $\tl$ is a rigid monoidal category.  On morphisms, the dual is constructed using $e_i$ and $c_i$ which has the effect of rotating diagrams by $180^\circ$.

As $\tl$ is a semisimple $\Cnum$-category (Theorem \ref{thm:tlss}), 
it is abelian.  We have just seen that $\tl$ is rigid, and by construction the monoidal functor on $\tl$ is bilinear on morphisms.  Finally, $\End_{\tl}(\ti)=\Cnum$, so $\tl$ is a semisimple tensor category, i.e., a fusion category, according to \cite[Definition 4.1.1]{egno}.

As $\tl$ is semisimple, 
it follows from Lemma \ref{lem:duals} that every object is self-dual, and so we get:
\begin{proposition}
$\tl$ is a pivotal category where $\iota_X:X\to X^{\vee\vee}=X$ is the identity map.
\end{proposition}

The following result will be useful in what follows.  A particular case of this result, when $k$ is a multiple of 4, appears in \cite[Lemma 2.4]{mps:D2n}. 
\begin{lemma} \label{lem:fkfk-identity}
The map $c_k\circ e_k:F_k\otimes F_k\to F_k\otimes F_k$ is equal to $(-1)^k \id_{F_k}\otimes \id_{F_k}$.
\end{lemma}
\begin{proof}
Since $F_k \otimes F_k \cong F_0$, which is simple, we have $\dim_\Cnum\Hom(F_k \otimes F_k,F_k \otimes F_k)=1$, thus 
\[c_k\circ e_k=\lambda\id_{F_k}\otimes \id_{F_k}\]
for some scalar $\lambda$. 
Now we cup and cap the last $k$ strings as follows:
\[(\id_{F_k}\otimes e_k )((c_k\circ e_k)\otimes \id_{F_k})(\id_{F_k}\otimes c_k) 
 = \lambda(\id_{F_k}\otimes e_k )(\id_{F_k}\otimes c_k) 
 \]
 The left hand side is the identity on $F_k$, by the duality equations, and by Proposition \ref{prop:jw-facts} the right hand side is $(-1)^k[k+1]\id_{F_k}\otimes\id_{F_k}$.  So as $[k+1]=1$ we have $\lambda=(-1)^k$.  
\end{proof}

\begin{lemma} \label{lem:rotate}
$(e_k\otimes \id_{F_{k-i}})(\id_{F_k}\otimes\pi_{i,k-i}\otimes\id_{F_k-i})(\id_{F_k}\otimes\id_{F_i}\otimes c_{k-i}) = (\id_{F_{k-i}}\otimes c_i)(\iota_{k-1,i}\otimes \id_{F_i})$
\end{lemma}
\begin{proof}
Use that $\id_{F_k}^*=\id_{F_k^*}=\id_{F_k}$ and the definition of a dual map using evaluation and coevaluation: see \cite[(2.47)]{egno}.
\end{proof}

We can illustrate Lemmas \ref{lem:fkfk-identity} and \ref{lem:rotate} as follows:
\[
\begin{tikzpicture}
[auto, block/.style ={rectangle, draw=black, thick, fill=blue!10, text width=2em, align=center, rounded corners, minimum height=2em},
 dot/.style = {circle,draw=black, thick, fill=green!10} ]
\node[block] (t1) at (0,2) {$F_{k}$};
\node[block] (t2) at (1.5,2) {$F_k$};
\node[block] (b1) at (0,0) {$F_k$};
\node[block] (b2) at (1.5,0) {$F_k$};
\draw[thick] (t1) to[out=-90,in=-90] (t2);
\draw[thick] (b1) to[out=90,in=90] (b2);
\node (s1) at (2.5,1) {$=$};
\node () at (3.4,1) {$(-1)^k$};
\node[block] (t3) at (4,2) {$F_{k}$};
\node[block] (t4) at (5.5,2) {$F_k$};
\node[block] (b3) at (4,0) {$F_k$};
\node[block] (b4) at (5.5,0) {$F_k$};
\draw[thick] (t3) to (b3);
\draw[thick] (t4) to (b4);
\node[block] (m1) at (9,1) {$F_k$};
\node (t5) at (10,2) {$k-i$};
\node (b5) at (8,0) {$k$};
\node (b6) at (8.8,0) {$i$};
\coordinate (i4) at (9.6,0.5) {};
\coordinate (i5) at (8.4,1.5) {};
\coordinate (i6) at (9,0.25) {};
\draw[thick] (m1) to[out=90,in=0] (i5) to[out=180,in=90] (b5);
\draw[thick, transform canvas={xshift=-1ex}] (m1) -- (i6);
\draw[thick, transform canvas={xshift=1ex}] (m1) to[out=-90,in=180] (i4) to[out=0,in=-90] (t5);
\node (s2) at (11,1) {$=$};
\node[block] (m2) at (12,1) {$F_k$};
\node (t8) at (12,2) {$k-i$};
\node (b8) at (12,0) {$k$};
\node (b9) at (12.8,0) {$i$};
\coordinate (i7) at (12.4,1.5) {};
\draw[thick] (m2) to[out=90,in=180] (i7) to[out=0,in=90] (b9);
\draw[thick] (m2) -- (b8);
\draw[thick, transform canvas={xshift=-1ex}] (m2) -- (t8);
\end{tikzpicture}
\]

We will also need the following structure.  Recall from Section \ref{ss:jw} that we have a fixed square root $t$ of $q$.  
\begin{proposition}[{\cite[Lemma XII.6.4.1]{turaev}}]\label{prop:braiding}
$\tl$ is a braided monoidal category with braiding $\sigma_{1,1}:1\otimes 1\to 1\otimes 1$ defined by 
$\sigma_{1,1}=t\id_{1\otimes 1}+t^{-1}U_{0,2}$.  
Its inverse is defined by 
$\sigma_{1,1}^{-1}=t^{-1}\id_{1\otimes 1}+tU_{0,2}$.
\end{proposition}
We braid other objects using the axioms for a braiding: for example, $\sigma_{1,2}=(\id_1\otimes\sigma_{1,1})\circ(\sigma_{1,1}\otimes\id_1)$.

We draw $\sigma$ going up the page, with the overcrossing going left to right:
\[
\begin{tikzpicture}
[auto,
 block/.style ={rectangle, draw=black, thick, fill=green!10, text width=2em, align=center, rounded corners, minimum height=2em},
kill/.style = {circle, 
 fill=white, minimum size=#1,
              inner sep=0pt, outer sep=0pt}
]
\node (t1) at (0,1) {$1$};
\node (t2) at (1,1) {$1$};
\node (b1) at (0,0) {$1$};
\node (b2) at (1,0) {$1$};
\node at (-0.8,0.5) {$\sigma_{1,1}\;=$};
\draw[thick] (b2) -- (t1);
\node[kill=7pt] (z3) at (0.5,0.5) {};
\draw[thick] (b1) -- (t2);
\node at (1.5,0.5) {$=$};
\node at (2.1,0.5) {$t$};
\node (t1) at (2.5,1) {$1$};
\node (t2) at (3,1) {$1$};
\node (b1) at (2.5,0) {$1$};
\node (b2) at (3,0) {$1$};
\draw[thick] (b1) -- (t1);
\draw[thick] (b2) -- (t2);
\node at (4.2,0.5) {$+\;\;\;\;t^{-1}$};
\node (t1) at (5,1) {$1$};
\node (t2) at (5.5,1) {$1$};
\node (b1) at (5,0) {$1$};
\node (b2) at (5.5,0) {$1$};
\draw[thick] (t1) to[out=-90,in=-90] (t2);
\draw[thick] (b1) to[out=90,in=90] (b2);
\end{tikzpicture}
\]

We will make use of the following result: see \cite[Corollary 4.6]{fy}.
\begin{theorem}[Freyd-Yetter]\label{thm:fy}
In a pivotal braided category, an equation of maps holds if and only if the tangles representing the maps are regularly isotopic.
\end{theorem}

\begin{lemma}\label{lem:crossJW}
We have 
$f_{a+b}\circ\sigma_{a,b}=t^{ab}f_{a+b}$ 
and 
$f_{a+b}\circ\sigma^{-1}_{a,b}=t^{-ab}f_{a+b}$ 
as maps $F_a\otimes F_b\to F_{a+b}$, and
$\sigma_{a,b}\circ f_{a+b}=t^{ab}f_{a+b}$ 
and
$\sigma^{-1}_{a,b}\circ f_{a+b}=t^{-ab}f_{a+b}$ 
as maps $F_{a+b}\to F_a\otimes F_b$.
\[
\begin{tikzpicture}
[auto, block/.style ={rectangle, draw=black, thick, fill=blue!10, text width=2em, align=center, rounded corners, minimum height=2em},
 wblock/.style ={rectangle, draw=black, thick, fill=blue!10, text width=2em, align=center, rounded corners, minimum height=2em, minimum width=4em},
kill/.style = {circle, 
 fill=white, minimum size=#1,
              inner sep=0pt, outer sep=0pt},
 dot/.style = {circle,draw=black, thick, fill=green!10} ]
\node (x1) at (0.3,1.15) {};
\node (y1) at (1.2,1.15) {};
\node[wblock] (p4) at (0.7,1.4) {$F_{a+b}$};
\node (a1) at (0.2,0) {$a$};
\node (b1) at (1.3,0) {$b$};
\draw[thick] (b1) -- (x1);
\node[kill=7pt] (z3) at (0.77,0.64) {};
\draw[thick] (a1) -- (y1);
\node at (2.7,0.7) {$=\;\;\;\;t^{ab}$};
\node (x2) at (3.8,1.15) {};
\node (y2) at (4.6,1.15) {};
\node[wblock] (p2) at (4.2,1.4) {$F_{a+b}$};
\node (a2) at (3.8,0) {$a$};
\node (b2) at (4.6,0) {$b$};
\draw[thick] (b2) -- (y2);
\draw[thick] (a2) -- (x2);
\end{tikzpicture}
\hspace{4em}
\begin{tikzpicture}
[auto, block/.style ={rectangle, draw=black, thick, fill=blue!10, text width=2em, align=center, rounded corners, minimum height=2em},
 wblock/.style ={rectangle, draw=black, thick, fill=blue!10, text width=2em, align=center, rounded corners, minimum height=2em, minimum width=4em},
kill/.style = {circle, 
 fill=white, minimum size=#1,
              inner sep=0pt, outer sep=0pt},
 dot/.style = {circle,draw=black, thick, fill=green!10} ]
\node (x1) at (0.3,0.35) {};
\node (y1) at (1.2,0.35) {};
\node[wblock] (p4) at (0.7,0.1) {$F_{a+b}$};
\node (a1) at (0.2,1.5) {$a$};
\node (b1) at (1.3,1.5) {$b$};
\draw[thick] (a1) -- (y1);
\node[kill=7pt] (z3) at (0.75,0.88) {};
\draw[thick] (b1) -- (x1);
\node at (2.7,0.7) {$=\;\;\;\;t^{ab}$};
\node (x2) at (3.8,0.25) {};
\node (y2) at (4.6,0.25) {};
\node[wblock] (p2) at (4.2,0) {$F_{a+b}$};
\node (a2) at (3.8,1.4) {$a$};
\node (b2) at (4.6,1.4) {$b$};
\draw[thick] (b2) -- (y2);
\draw[thick] (a2) -- (x2);
\end{tikzpicture}
\]
\end{lemma}
\begin{proof}
These follow immediately from the definition of the braiding in Proposition \ref{prop:braiding} and parts (ii) and (iii) of Proposition \ref{prop:jw-facts}.
\end{proof}

In a braided pivotal category we have left and right twists $\theta^\ell_X,\theta^r_X:X\to X$ for each object $X$, defined by:
\[ \theta^\ell_X=(e_X\otimes\id_X)\circ (\iota_{^\vee X}\otimes\sigma_{X,X})\circ(c^r_{X}\otimes\id_X) \;\;\text{ and } \;\;
\theta^r_X=(\id_X\otimes e^r_{X})\circ (\sigma_{X,X}\otimes\iota_{^\vee X}^{-1})\circ(\id_X\otimes c_X).\]
There are explicit descriptions for their inverses and duals: see \cite[Section 3.3.1]{tv}.

The twists for $\tl$ were studied in \cite[Section XII.6]{turaev}.  We calculate their effect on the Jones-Wenzl idempotents.
\begin{lemma}\label{lem:twistJW}
In $\tl$ we have $\theta^\ell_{F_j}=\theta^r_{F_j}=(-1)^jt^{j^2+2j}\id_{F_j}$.
\end{lemma}
Graphically, we have:
\[
\begin{tikzpicture}
[auto, block/.style ={rectangle, draw=black, thick, fill=blue!10, text width=2em, align=center, rounded corners, minimum height=2em},
 wblock/.style ={rectangle, draw=black, thick, fill=blue!10, text width=2em, align=center, rounded corners, minimum height=2em, minimum width=4em},
kill/.style = {circle, 
 fill=white, minimum size=#1,
              inner sep=0pt, outer sep=0pt},
 dot/.style = {circle,draw=black, thick, fill=green!10} ]
\node[block] (bo) at (0,0.5) {$F_j$};
\node[block] (t) at (0,2.5) {$F_j$};
\coordinate (b) at (0.3,1.5) {};
\coordinate (c) at (0.9,1.5) {};
\coordinate (c') at (0.9,1.5) {};

\draw[thick] (c') to[out=90,in=50] (b) to[out=-130,in=90] (bo); 
\node[kill=7pt] at (b) {};
\draw[thick] 
(t) to[out=-90,in=130] (b) to[out=-50,in=-90] (c);
\node at (2.2,1.5) {$=(-1)^jt^{-j^2-2j}$}; 
\node[block] (bo2) at (3.5,0.5) {$F_j$};
\node[block] (t2) at (3.5,2.5) {$F_j$};
\draw[thick] (bo2) to (t2);

\node[block] (bo) at (6,0.5) {$F_j$};
\node[block] (t) at (6,2.5) {$F_j$};
\coordinate (b) at (6.3,1.5) {};
\coordinate (c') at (6.9,1.5) {};
\coordinate (c) at (6.9,1.5) {};

\draw[thick] 
(t) to[out=-90,in=130] (b) to[out=-50,in=-90] (c);
\node[kill=7pt] at (b) {};
\draw[thick] (c') to[out=90,in=50] (b) to[out=-130,in=90] (bo); 
\node at (8.1,1.5) {$=(-1)^jt^{j^2+2j}$}; 
\node[block] (bo2) at (9.2,0.5) {$F_j$};
\node[block] (t2) at (9.2,2.5) {$F_j$};
\draw[thick] (bo2) to (t2);

\end{tikzpicture}
\]
\begin{proof}
We prove the statement for $(\theta^\ell_{F_j})^{-1}$; the others are similar, or follow from \cite[Lemma 3.2]{tv}.
We use induction on $j$.  The base case follows directly from the definition of the braiding in Proposition \ref{prop:braiding}:
\[ (\theta^\ell_{F_1})^{-1}=(t^{-1}\delta+t)\id_{F_1}=-t^-t^{-3}+t=-t^{-3}.\]
Then the inductive step goes as follows:
\[
\begin{tikzpicture}
[auto, block/.style ={rectangle, draw=black, thick, fill=blue!10, text width=2em, align=center, rounded corners, minimum height=2em},
 wblock/.style ={rectangle, draw=black, thick, fill=blue!10, text width=2em, align=center, rounded corners, minimum height=2em, minimum width=4em},
kill/.style = {circle, 
 fill=white, minimum size=#1,
              inner sep=0pt, outer sep=0pt},
 dot/.style = {circle,draw=black, thick, fill=green!10} ]
\node[block] (b0) at (0,0) {$F_{j+1}$};
\node[block] (t0) at (0,3.5) {$F_{j+1}$};

\node (tj) at (-0.2,3.29) {};
\node (t1) at (0.2,3.29) {};
\node (bj) at (-0.2,0.21) {};
\node (b1) at (0.2,0.21) {};

\coordinate (a') at (0.8,0.8) {};
\coordinate (a) at (1,1.2) {};
\coordinate (c) at (1.2,2.1) {};
\coordinate (c') at (1,2.6) {};

\coordinate (b) at (1.5,1.75) {};
\coordinate (b') at (2,1.75) {};
\coordinate (bb) at (1.5,1.75) {};
\coordinate (bb') at (2,1.75) {};

\node  at (-0.4,2.9) {$j$};
\node  at (0.4,2.9) {$1$};

\draw[thick] (bb') to [out=90,in=0] (c') to[out=180,in=90] (bj);
\draw[thick] (bb) to [out=90,in=0] (c) to[out=180,in=90] (b1);

\node[kill=8pt] at (0.2,2.18) {};
\node[kill=8pt] at (-0.04,1.6) {};
\node[kill=8pt] at (0.37,1.6) {};
\node[kill=8pt] at (0.26,1.05) {};

\draw[thick] (tj) to[out=-90,in=180] (a') to[out=0,in=-90] (b');
\draw[thick] (t1) to[out=-90,in=180] (a) to[out=0,in=-90] (b); 

\node at (3,1.75) {$= \;\; -t^{-3}$};

\node[block] (b0) at (4,0) {$F_{j+1}$};
\node[block] (t0) at (4,3.5) {$F_{j+1}$};

\node (tj) at (3.8,3.29) {};
\node (t1) at (4.2,3.29) {};
\node (bj) at (3.8,0.21) {};
\node (b1) at (4.2,0.21) {};

\coordinate (a') at (4.8,0.8) {};
\coordinate (c') at (4.8,2.7) {};

\coordinate (b') at (5.5,1.75) {};
\coordinate (bb') at (5.5,1.75) {};

\node  at (3.6,2.9) {$j$};
\node  at (4.4,2.9) {$1$};

\draw[thick] (bb') to [out=90,in=0] (c') to[out=180,in=90] (bj);
\draw[thick] (4.37,1.5) to[out=-90,in=90] (b1);

\node[kill=8pt] at (4.26,2.5) {};
\node[kill=10pt] at (3.9,1.7) {};
\node[kill=8pt] at (4.3,0.92) {};

\draw[thick] (tj) to[out=-90,in=180] (a') to[out=0,in=-90] (b');
\draw[thick] (t1) to[out=-90,in=90] (4.37,1.5);

\node at (6.6,1.75) {$= \;\; -t^{-3}$};

\node[block] (b0) at (8,0) {$F_{j+1}$};
\node[block] (t0) at (8,3.5) {$F_{j+1}$};

\node (tj) at (7.8,3.29) {};
\node (t1) at (8.2,3.29) {};
\node (bj) at (7.8,0.21) {};
\node (b1) at (8.2,0.21) {};

\coordinate (a') at (8.7,1.2) {};
\coordinate (c') at (8.7,2.3) {};

\coordinate (b') at (9.2,1.75) {};
\coordinate (bb') at (9.2,1.75) {};

\node  at (7.6,2.9) {$j$};
\node  at (8.4,2.9) {$1$};

\draw[thick] (bb') to [out=90,in=0] (c') to[out=180,in=90] (bj);
\node[kill=7pt] at (7.84,0.92) {};
\draw[thick] (7.5,1.5) to[out=-90,in=90] (b1);

\node[kill=7pt] at (7.980,1.75) {};

\draw[thick] (t1) to[out=-90,in=90] (7.5,1.5); 
\node[kill=7pt] at (7.83,2.28) {};
\draw[thick] (tj) to[out=-90,in=180] (a') to[out=0,in=-90] (b');

\node at (10.96,1.75) {$= \;\; (-1)^{j+1}t^{-3-j^2-2j}$};

\node[block] (b0) at (13,0) {$F_{j+1}$};
\node[block] (t0) at (13,3.5) {$F_{j+1}$};

\node (tj) at (12.8,3.29) {};
\node (t1) at (13.2,3.29) {};
\node (bj) at (12.8,0.21) {};
\node (b1) at (13.2,0.21) {};

\node (l) at (12.8,1.62) {};
\node (r) at (13.2,1.62) {};
\node (l') at (12.8,1.88) {};
\node (r') at (13.2,1.88) {};

\node  at (12.6,2.9) {$j$};
\node  at (13.4,2.9) {$1$};

\draw[thick] (r') to [out=-90,in=90] (bj);
\node[kill=7pt] at (13,1.04) {};
\draw[thick] (l') to[out=-90,in=90] (b1);

\draw[thick] (t1) to[out=-90,in=90] (l); 
\node[kill=7pt] at (13,2.46) {};
\draw[thick] (tj) to[out=-90,in=90] (r);

\end{tikzpicture}
\]
We remove the twists using Lemma \ref{lem:crossJW} to get an additional factor of $t^{-2j}$, giving overall scalar $(-1)^{j+1}t^{-j^2-4j-3}=(-1)^{j+1}t^{-(j+1)^2-2(j+1)}$.
\end{proof}

Turaev showed \cite[Theorem XII.6.6]{turaev} that $\tl$ is a ribbon category, following the definition in \cite[Section I.1.4]{turaev}.  Note that Lemma \ref{lem:twistJW} implies $\theta^\ell=\theta^r$, so we see directly that $\tl$ satisfies the definition in \cite[Section 3.3.2]{tv}.

\subsection{Twisted Frobenius structure}

We have $k$ fixed for $\tl=\tl_k$.  Throughout this section we will work with natural numbers $i$ and $j$ such that $i+j=k$.
\begin{definition}
Let $W=F_k$.  Let $n:\Sigma\to F_k$ be defined by projection.
\end{definition}
\begin{lemma}\label{lem:sigma-frob}
$(F_k,n)$ is a Frobenius structure on $\Sigma$.
\end{lemma}
\begin{proof}
We need to show that the map
\[ \varphi:\Sigma \larr{\id_\Sigma\otimes c_\Sigma} \Sigma\otimes \Sigma\otimes \Sigma^\vee\larr{m\otimes \id_{\Sigma^\vee}}\Sigma\otimes \Sigma^\vee\larr{n\otimes \id_{\Sigma^\vee}} F_k\otimes \Sigma^\vee \]
is invertible.  We split the map into components
\[ \varphi_{i,j}:F_i \larr{\id_{F_i}\otimes c_{j}} F_i\otimes F_j\otimes F_j\larr{m\otimes \id_{F_j}}F_{i+j}\otimes F_j\larr{n\otimes \id_{F_j}} F_k\otimes F_j \]
and notice that the final map $n\otimes \id_{F_j}$ is nonzero iff $i+j=k$, in which case $n=\id_{F_k}:F_k\to F_k$.  Then $m=\pi_{i,k-i}$.  So the nonzero components are
\[ \varphi_i:=\varphi_{i,k-i}:F_i \larr{\id_{F_i}\otimes c_{k-i}} F_i\otimes F_{k-i}\otimes F_{k-i}\larr{\pi_{i,k-i}\otimes \id_{F_{k-i}}}F_k\otimes F_{k-i}. \]
We define a map in the opposite direction with components as follows:
\[ \frac{(-1)^{k-i}}{[i+1]}\psi_i:=
F_k\otimes F_{k-i}\arrr{\iota_{i,k-1}\otimes\id_{F_{k-i}}}F_i\otimes F_{k-i}\otimes F_{k-i}\arrr{\id_{F_i}\otimes e_{k-1}}F_i. \]
Graphically, we have:
\[
\begin{tikzpicture}
[auto, block/.style ={rectangle, draw=black, thick, fill=blue!10, text width=2em, align=center, rounded corners, minimum height=2em},
 dot/.style = {circle,draw=black, thick, fill=green!10} ]
\node[block] (b1) at (0.5,0) {$F_i$};
\node[block] (t1) at (0.5,1.5) {$F_k$};
\node[block] (t2) at (1.8,1.5) {$F_{k-i}$};
\node (i2) at (1.6,1.28) {};
\node (ca) at (1,1.18) {};
\node (sa) at (-0.5,0.7) {$\varphi_i =$};
\draw[thick, transform canvas={xshift=-1ex}] (b1) -- (t1);
\draw[thick, transform canvas={xshift=1ex}] (t1) to[out=-90,in=-90] (i2);
\node (sb) at (4.7,0.7) {$\psi_i =(-1)^{k-i}[i+1]$};
\node[block] (b1) at (6.5,1.5) {$F_i$};
\node[block] (t1) at (6.5,0) {$F_k$};
\node[block] (t2) at (7.8,0) {$F_{k-i}$};
\node (i2) at (7.6,0.22) {};
\node (ca) at (7,0.32) {};
\draw[thick, transform canvas={xshift=-1ex}] (b1) -- (t1);
\draw[thick, transform canvas={xshift=1ex}] (t1) to[out=90,in=90] (i2);
\end{tikzpicture}
\]
Then $\psi_i\circ\varphi_i=\id_{F_i}$ by Proposition \ref{prop:jw-facts}(v), using $[j]=[k+2-j]$ and $[1]=1$.

To compute $\varphi_i\circ\psi_i$ we use Lemmas \ref{lem:fkfk-identity} and \ref{lem:rotate} as follows:
\[
\begin{tikzpicture}
[auto, block/.style ={rectangle, draw=black, thick, fill=blue!10, text width=2em, align=center, rounded corners, minimum height=2em},
 dot/.style = {circle,draw=black, thick, fill=green!10} ]
\node[block] (e1) at (0,3.2) {$F_k$};
\node[block] (e2) at (1.2,3.2) {$F_{k-i}$};
\node[block] (a1) at (0,0) {$F_k$};
\node[block] (a2) at (1.2,0) {$F_{k-i}$};
\draw[thick, transform canvas={xshift=-1ex}] (a1) -- (e1);
\draw[thick, transform canvas={xshift=1ex}] (e1) to[out=-90,in=-90] (e2);
\draw[thick, transform canvas={xshift=1ex}] (a1) to[out=90,in=90] (a2);
\node[block] (e3) at (3,3.2) {$F_k$};
\node[block] (e4) at (5.5,3.2) {$F_{k-i}$};
\node[block] (d3) at (4.5,2.2) {$F_k$};
\node[block] (b3) at (4.5,1) {$F_k$};
\node[block] (a3) at (3,0) {$F_k$};
\node[block] (a4) at (5.5,0) {$F_{k-i}$};
\coordinate (i3) at (3.8,0.5) {};
\coordinate (i4) at (3.8,2.7) {};
\draw[thick, transform canvas={xshift=-1ex}] (b3) -- (d3);
\draw[thick] (d3) to[out=-45,in=-90] (e4);
\draw[thick] (b3) to[out=45,in=90] (a4);
\draw[thick] (e3) to[out=-90,in=-90] (i4) to[out=90,in=90] (d3);
\draw[thick] (a3) to[out=90,in=90] (i3) to[out=-90,in=-90] (b3);
\node[block] (e5) at (7.3,3.2) {$F_k$};
\node[block] (e6) at (9.5,3.2) {$F_{k-i}$};
\node[block] (d5) at (8.5,2.2) {$F_k$};
\node[block] (b5) at (8.5,1) {$F_k$};
\node[block] (a5) at (7.3,0) {$F_k$};
\node[block] (a6) at (9.5,0) {$F_{k-i}$};
\coordinate (i5) at (7.8,0.8) {};
\coordinate (i6) at (7.8,2.7) {};
\draw[thick, transform canvas={xshift=-1ex}] (b5) -- (d5);
\draw[thick] (d5) to[out=-45,in=-90] (e6);
\draw[thick] (b5) to[out=45,in=90] (a6);
\draw[thick] (e5) to (a5);
\draw[thick] (b5) to[out=-90,in=-90] (i5) to [out=90,in=-90] (i6) to[out=90,in=90] (d5);
\node[block] (e7) at (11.3,3.2) {$F_k$};
\node[block] (e8) at (12.5,3.2) {$F_{k-i}$};
\node[block] (d7) at (12.5,2.2) {$F_k$};
\node[block] (b7) at (12.5,1) {$F_k$};
\node[block] (a7) at (11.3,0) {$F_k$};
\node[block] (a8) at (12.5,0) {$F_{k-i}$};
\coordinate (i7) at (13.3,1) {};
\coordinate (i8) at (13.3,2.4) {};
\draw[thick] (b7) -- (d7);
\draw[thick] (e7) to (a7);
\draw[thick, transform canvas={xshift=-1ex}] (e8) to (d7);
\draw[thick, transform canvas={xshift=-1ex}] (a8) to (b7);
\draw[thick] (d5) to[out=-45,in=-90] (e6);
\draw[thick] (b5) to[out=45,in=90] (a6);
\draw[thick] (e7) to (a7);
\draw[thick] (b7) to[out=-45,in=-90] (i7) to [out=90,in=-90] (i8) to[out=90,in=45] (d7);
\node at (2,1.6) {$=$};
\node at (6.3,1.6) {$=$};
\node at (10.3,1.6) {$=$};
\node at (0,1.6) {$i$};
\node at (4.5,1.6) {$i$};
\node at (8.5,1.6) {$i$};
\node at (13.5,1.6) {$i$};
\end{tikzpicture}
\]
Then using Proposition \ref{prop:jw-facts}(v) again, we get $\varphi_i\circ\psi_i=\id_{F_k}\otimes \id_{F_{k-1}}$.  So $\varphi$ and $\psi$ are inverse, and therefore $\varphi$ is invertible.
\end{proof}

Recall that dualizing a map in $TL$ rotates by $180^\circ$.   So we can construct the $i$th graded part of the Nakayama morphism $\alpha$ of $\Sigma$ directly from Definition \ref{def:nak} as follows:
\[
\begin{tikzpicture}
[auto, block/.style ={rectangle, draw=black, thick, fill=blue!10, text width=2em, align=center, rounded corners, minimum height=2em},
 dot/.style = {circle,draw=black, thick, fill=green!10} ]
\node[block] (b1) at (0,0) {$F_k$};
\node[block] (b2) at (1.5,0) {$F_{i}$};
\node[block] (m1) at (1.5,1.5) {$F_k$};
\node[block] (m2) at (3,1.5) {$F_{k-i}$};
\node[block] (t1) at (1.5,3) {$F_i$};
\node[block] (t2) at (3,3) {$F_k$};
\draw[thick, transform canvas={xshift=-1ex}] (b2) -- (m1);
\draw[thick, transform canvas={xshift=1ex}] (m2) -- (t2);
\draw[thick, transform canvas={xshift=1ex}] (b1) to[out=90,in=135] (m1);
\draw[thick, transform canvas={xshift=1ex}] (t1) to[out=-70,in=-110] (t2);
\draw[thick, transform canvas={xshift=1ex}] (m1) to[out=-70,in=-110] (m2);
\node at (-1.2,2) {$\alpha\;=\;(-1)^{k-i}[i+1]$};
\end{tikzpicture}
\]
\begin{theorem}\label{thm:alpha^2}
$\alpha^2=1$.
\end{theorem}
\begin{proof}
Following Definition \ref{def:nak-order}, we should calculate the composition
\[A=\ti\otimes A\to V^{\otimes 2}\otimes A\to A\otimes V^{\otimes 2}\to A\otimes\ti=A.\] 
Let $j=k-i$.  Omitting the scalar $[i+1]^2$, we get:
\[
\begin{tikzpicture}
[auto, block/.style ={rectangle, draw=black, thick, fill=blue!10, text width=2em, align=center, rounded corners, minimum height=2em},
 wblock/.style ={rectangle, draw=black, thick, fill=blue!10, text width=2em, align=center, rounded corners, minimum height=2em, minimum width=4em},
kill/.style = {circle, 
 fill=white, minimum size=#1,
              inner sep=0pt, outer sep=0pt},
 dot/.style = {circle,draw=black, thick, fill=green!10} ]
\node[block] (a0) at (-1.5,0) {$F_k$};
\node[block] (a1) at (0,0) {$F_k$};
\node[block] (a2) at (1.5,-0.7) {$F_{i}$};
\node[block] (b1) at (1.5,1.5) {$F_k$};
\node[block] (b2) at (3,1.5) {$F_{j}$};
\coordinate (bc) at (1.5,2.3) {};
\node[block] (c1) at (0,3.2) {$F_i$};
\node[block] (c2) at (3,3.2) {$F_k$};
\node[block] (d1) at (0,4.5) {$F_k$};
\node[block] (d2) at (1.5,4.5) {$F_{j}$};
\node[block] (e1) at (0,6.7) {$F_i$};
\node[block] (e2) at (1.5,6) {$F_{k}$};
\draw[thick, transform canvas={xshift=-1ex}] (a2) -- (b1);
\draw[thick, transform canvas={xshift=1ex}] (b2) -- (c2);
\draw[thick] (a1) to[out=90,in=130] (b1);
\draw[thick] (c1) to[out=-70,in=180] (bc) to[out=0,in=-110] (c2);
\draw[thick, transform canvas={xshift=1ex}] (b1) to[out=-70,in=-110] (b2);
\draw[thick, transform canvas={xshift=-1ex}] (c1) -- (d1);
\draw[thick, transform canvas={xshift=1ex}] (d1) to[out=-70,in=-110] (d2);
\draw[thick, transform canvas={xshift=-1ex}] (e1) to[out=-70,in=-110] (e2);
\draw[thick, transform canvas={xshift=1ex}] (d2) -- (e2);
\draw[thick] (a0) to[out=90,in=130] (d1);
\draw[thick] (a0) to[out=-70,in=-110] (a1);
\draw[thick] (e2) to[out=60,in=90] (c2);
\node at (4.3,3) {$=$};
\node[block] (f) at (5,-0.7) {$F_i$};
\node[wblock] (g) at (5.5,1.5) {$F_k$};
\node[wblock] (i) at (6.5,4.5) {$F_k$};
\node[block] (j) at (7,6.7) {$F_{i}$};
\draw[thick, transform canvas={xshift=-1ex}] (f) -- (g);
\draw[thick, transform canvas={xshift=1ex}] (i) -- (j);
\draw[thick] (g) -- (i);
\draw[thick, transform canvas={xshift=-2ex}] (g) to[out=100,in=100] (i);
\draw[thick, transform canvas={xshift=2ex}] (g) to[out=-80,in=-80] (i);
\node at (6.2,3) {$i$};
\node at (6.9,2.2) {$j$};
\node at (5.55,3.8) {$j$};
\end{tikzpicture}
\]
where we have used the absorption property and the duality relations, and labelled the curves by the number of strands.
Now, using more duality relations, this is equal to the first of the following diagrams:
\[
\begin{tikzpicture}
[auto, block/.style ={rectangle, draw=black, thick, fill=blue!10, text width=2em, align=center, rounded corners, minimum height=2em},
 wblock/.style ={rectangle, draw=black, thick, fill=blue!10, text width=2em, align=center, rounded corners, minimum height=2em, minimum width=4em},
kill/.style = {circle, 
 fill=white, minimum size=#1,
              inner sep=0pt, outer sep=0pt},
 dot/.style = {circle,draw=black, thick, fill=green!10} ]
\node[block] (f2) at (2,0) {$F_i$};
\node[wblock] (g2) at (2.2,1.5) {$F_k$};
\node[wblock] (i2) at (0,5.5) {$F_k$};
\node[block] (j2) at (0.2,7) {$F_{i}$};
\draw[thick, transform canvas={xshift=-1ex}] (f2) to (g2);
\draw[thick, transform canvas={xshift=1ex}] (i2) to (j2);
\coordinate (a) at (-1,6.2) {};
\coordinate (b) at (-1.8,3.5) {};
\coordinate (c) at (-1,0.8) {};
\coordinate (d) at (0.5,2.2) {};
\draw[thick, transform canvas={xshift=-1ex}] (i2) to[out=110,in=0] (a) to[out=180,in=90] (b) to[out=-90,in=180] (c) to [out=0,in=180] (d) to[out=-10,in=90] (g2);
\coordinate (e) at (3.2,0.8) {};
\coordinate (f) at (4,3.5) {};
\coordinate (g) at (3.2,6.2) {};
\coordinate (h) at (1.7,4.8) {};
\draw[thick, transform canvas={xshift=+1ex}] (g2) to[out=-70,in=180] (e) to[out=0,in=-90] (f) to[out=90,in=0] (g) to [out=180,in=0] (h) to[out=170,in=-90] (i2);
\coordinate (z') at (1,2.5) {};
\coordinate (a') at (-0.8,2) {};
\coordinate (b') at (-1.5,2.2) {};
\coordinate (c') at (0,3.5) {};
\coordinate (d') at (2.6,3.5) {};
\coordinate (e') at (3.55,1.5) {};
\coordinate (e'') at (3.3,1.5) {};
\coordinate (f') at (4,3.5) {};
\coordinate (g') at (3.5,5.7) {};
\coordinate (h') at (2,4.5) {};
\draw[thick, transform canvas={xshift=+1ex}] (g2) to[out=80,in=0] (z') to[out=180,in=60] (a') to[out=-120,in=-90] (b') to[out=90,in=180] (c') to[out=0,in=180] (d') to[out=0,in=180] (e'');
\draw[thick, transform canvas={xshift=-1ex}] (e') to[out=0,in=-90] (f') to[out=90,in=0] (g') to [out=180,in=0] (h') to[out=170,in=-90] (i2);
\node at (-1.6,5.6) {$j$};
\node at (-0.7,3) {$i$};
\node at (4.1,1) {$j$};
\node at (4.8,3.5) {$=$};
\node[block] (f2) at (9.5,0) {$F_i$};
\node[wblock] (g2) at (9.7,1.5) {$F_k$};
\node[wblock] (i2) at (7.5,5.5) {$F_k$};
\node[block] (j2) at (7.7,7) {$F_{i}$};
\coordinate (a) at (6.5,6.2) {};
\coordinate (b) at (5.7,3.5) {};
\coordinate (c) at (6.5,0.8) {};
\coordinate (d) at (8,2.2) {};
\coordinate (e) at (10.7,0.8) {};
\coordinate (f) at (11.5,3.5) {};
\coordinate (g) at (10.7,6.2) {};
\coordinate (h) at (9.2,4.8) {};
\coordinate (z') at (8.5,2.5) {};
\coordinate (a') at (6.3,1.6) {};
\coordinate (ab') at (5.8,1.9) {};
\coordinate (b') at (6.4,2.3) {};
\coordinate (b'') at (6.1,2.3) {};
\coordinate (c') at (7.8,0.5) {};
\coordinate (d') at (10.8,1.1) {};
\coordinate (e') at (11.3,1.5) {};
\coordinate (f') at (11.5,3.5) {};
\coordinate (g') at (11,5.7) {};
\coordinate (h') at (9.5,4.5) {};
\draw[thick, transform canvas={xshift=-1ex}] (b') to[out=0,in=150] (c') to[out=-30,in=200] (d') to[out=20,in=-120] (e') to[out=60,in=-90] (f') to[out=90,in=0] (g') to [out=180,in=0] (h') to[out=170,in=-90] (i2);
\node[kill=7pt] at (7,1.55) {};
\node[kill=7pt] at (9.42,0.6) {};
\node[kill=7pt] at (10.15,0.92) {};
\node[kill=7pt] at (6.8,2) {};

\draw[thick, transform canvas={xshift=-1ex}] (i2) to[out=110,in=0] (a) to[out=180,in=90] (b) to[out=-90,in=180] (c) to [out=0,in=180] (d) to[out=-10,in=90] (g2);
\draw[thick, transform canvas={xshift=+1ex}] (g2) to[out=-70,in=180] (e) to[out=0,in=-90] (f) to[out=90,in=0] (g) to [out=180,in=0] (h) to[out=170,in=-90] (i2);
\draw[thick, transform canvas={xshift=-1ex}] (f2) to (g2);
\draw[thick, transform canvas={xshift=1ex}] (i2) to (j2);

\draw[thick, transform canvas={xshift=+1ex}] (g2) to[out=80,in=0] (z') to[out=180,in=60] (a') to[out=-120,in=-90] (ab') to [out=90,in=180] (b'');
\node at (6,5.5) {$j$};
\node at (6.3,2) {$i$};
\node at (11.5,1) {$j$};

\end{tikzpicture}
\]
Now we use the braiding, and appeal to Theorem \ref{thm:fy} to show the first diagram above is equal to the second.

From here we use Lemma \ref{lem:fkfk-identity} to obtain the first of the following diagrams, and Theorem \ref{thm:fy} again to obtain the second:
\[
\begin{tikzpicture}
[auto, block/.style ={rectangle, draw=black, thick, fill=blue!10, text width=2em, align=center, rounded corners, minimum height=2em},
 wblock/.style ={rectangle, draw=black, thick, fill=blue!10, text width=2em, align=center, rounded corners, minimum height=2em, minimum width=4em},
kill/.style = {circle, 
 fill=white, minimum size=#1,
              inner sep=0pt, outer sep=0pt},
 dot/.style = {circle,draw=black, thick, fill=green!10} ]
\node[block] (f2) at (9.5,0) {$F_i$};
\node[wblock] (g2) at (9.7,1.5) {$F_k$};
\node[wblock] (i2) at (7.5,5.5) {$F_k$};
\node (i2') at (8,5.28) {$$};
\node (g2') at (9.3,1.72) {$$};
\node[block] (j2) at (7.7,7) {$F_{i}$};
\coordinate (a) at (6.5,6.2) {};
\coordinate (b) at (5.7,3.5) {};
\coordinate (c) at (6.5,0.8) {};
\coordinate (d) at (8,2.2) {};
\coordinate (e) at (10.9,0.5) {};
\coordinate (f) at (12,3.5) {};
\coordinate (g) at (10.7,6.2) {};
\coordinate (h) at (9.2,4.8) {};
\coordinate (a') at (6.4,1.5) {};
\coordinate (ab') at (6,1.95) {};
\coordinate (b') at (6.15,2.3) {};
\coordinate (b'') at (6.5,2.3) {};
\coordinate (c') at (8.2,1) {};
\coordinate (e') at (10.5,0.8) {};
\coordinate (f') at (11.65,3.55) {};
\coordinate (g') at (11,5.7) {};
\coordinate (h') at (9.5,4.5) {};
\draw[thick, transform canvas={xshift=+1ex}] (g2) to[out=90,in=180] (g') to[out=0,in=90] (f') to[out=-90,in=0] (e') to[out=180,in=-10] (c') to[out=170,in=0] (b');
\node[kill=7pt] at (7.58,1.55) {};
\node[kill=7pt] at (9.45,0.85) {};
\node[kill=7pt] at (10.16,0.82) {};
\node[kill=7pt] at (7.1,1.9) {};
\draw[thick, transform canvas={xshift=-1ex}] (i2) to[out=110,in=0] (a) to[out=180,in=90] (b) to[out=-90,in=180] (c) to [out=0,in=-90] (i2');
\draw[thick, transform canvas={xshift=+1ex}] (g2) to[out=-70,in=180] (e)  to[out=0,in=-90] (f) to[out=90,in=0] (g) to [out=180,in=90] (g2');
\draw[thick, transform canvas={xshift=-1ex}] (f2) to (g2);
\draw[thick, transform canvas={xshift=1ex}] (i2) to (j2);
\draw[thick, transform canvas={xshift=-1ex}] (i2) to[out=-90,in=0] (a') to[out=-180,in=-90] (ab') to [out=90,in=180] (b'');
\node at (6,5.5) {$j$};
\node at (6.3,2) {$i$};
\node at (12.1,1) {$j$};
\node at (13.4,3.5) {$=$};
\node[block] (f2) at (16.5,0) {$F_i$};
\node[wblock] (g2) at (16.7,1.5) {$F_k$};
\node[wblock] (i2) at (16,5.5) {$F_k$};
\node (i2') at (16.5,5.28) {$$};
\node (g2') at (16.3,1.72) {$$};
\node[block] (j2) at (16.2,7) {$F_{i}$};
\coordinate (a) at (15.2,6.2) {};
\coordinate (b) at (14.7,5.5) {};
\coordinate (b0) at (14.7,5.5) {};
\coordinate (c) at (15.4,4.7) {};
\coordinate (d) at (17,2.2) {};
\coordinate (e) at (17.6,0.8) {};
\coordinate (f) at (18,2) {};
\coordinate (g) at (17,2.6) {};
\coordinate (a') at (15.3,3) {};
\coordinate (ab') at (14.8,3.5) {};
\coordinate (b') at (15,3.8) {};
\coordinate (b'') at (15.6,3.8) {};
\coordinate (e') at (17.3,4.1) {};
\coordinate (f') at (18,4.45) {};
\coordinate (f'') at (18,4.45) {};
\coordinate (g') at (17.6,4.8) {};

\draw[thick, transform canvas={xshift=-1ex}] (i2) to[out=110,in=0] (a) to[out=180,in=90] (b);

\draw[thick, transform canvas={xshift=+1ex}]  (f'') to[out=-90,in=0] (e') to[out=180,in=00] (b');

\node[kill=7pt] at (15.85,3.85) {};

\draw[thick, transform canvas={xshift=-1ex}] (i2) to[out=-90,in=0] (a') to[out=-180,in=-90] (ab') to [out=90,in=180] (b'');

\node[kill=7pt] at (15.85,4.72) {};

\draw[thick, transform canvas={xshift=-1ex}] (b0) to[out=-90,in=180] (c) to [out=0,in=-90] (i2');

\node[kill=7pt] at (16.9,4) {};

\draw[thick, transform canvas={xshift=+1ex}] (g2) to[out=90,in=180] (g') to[out=0,in=90] (f');

\node[kill=7pt] at (16.8,2.47) {};

\draw[thick, transform canvas={xshift=+1ex}] (g2) to[out=-70,in=180] (e)  to[out=0,in=-90] (f) to[out=90,in=0] (g) to [out=180,in=90] (g2');

\draw[thick, transform canvas={xshift=-1ex}] (f2) to (g2);
\draw[thick, transform canvas={xshift=1ex}] (i2) to (j2);
\node at (14.8,5.7) {$j$};
\node at (15,3.4) {$i$};
\node at (18.25,1) {$j$};

\end{tikzpicture}
\]
Now we use Lemma \ref{lem:crossJW} twice to remove the cossings next to the $F_k$ boxes, which introduces a scalar $t^{2ij}$.  Then we use Proposition \ref{prop:jw-facts} to remove the $F_k$ boxes, which introduces a scalar $1/[i+1]^2$, but this cancels with the scalar $[i+1]^2$ we omitted at the start of our calculation.  So we have so far shown that $\alpha^2$ is the following map:
\[
\begin{tikzpicture}
[auto, block/.style ={rectangle, draw=black, thick, fill=blue!10, text width=2em, align=center, rounded corners, minimum height=2em},
 wblock/.style ={rectangle, draw=black, thick, fill=blue!10, text width=2em, align=center, rounded corners, minimum height=2em, minimum width=4em},
kill/.style = {circle, 
 fill=white, minimum size=#1,
              inner sep=0pt, outer sep=0pt},
 dot/.style = {circle,draw=black, thick, fill=green!10} ]
\node[block] (bo) at (1,0) {$F_i$};
\node[block] (t) at (0,3) {$F_i$};
\coordinate (a) at (0,1.5) {};
\coordinate (b) at (-0.4,0.7) {};
\coordinate (c) at (-0.7,1.1) {};
\coordinate (c') at (-0.7,1.1) {};
\coordinate (d) at (0.5,1.5) {};
\coordinate (e) at (1,1.5) {};
\coordinate (f) at (1.7,1.9) {};
\coordinate (f') at (1.7,1.9) {};
\coordinate (g) at (1.4,2.3) {};

\draw[thick] (c') to[out=90,in=180] (a) to[out=0,in=180] (d) to[out=0,in=180] (e) to[out=0,in=-90] (f');
\node[kill=7pt] at (a) {};
\node[kill=7pt] at (e) {};
\draw[thick] (t) to[out=-90,in=90] (a) to[out=-90,in=0] (b) to[out=180,in=-90] (c);
\draw[thick] (bo) to[out=90,in=-90] (e) to[out=90,in=180] (g) to[out=0,in=90] (f);
\node at (-0.6,1.9) {$t^{2ij}$};
\node at (-2.2,1.5) {$\alpha^2\;\;\;\;=$};

\end{tikzpicture}
\]
Then by Lemma \ref{lem:twistJW} we can remove the twists at a cost of $t^{2i^2+4i}$, and we get $\alpha_i^2=\lambda_i\id_{F_i}$ where
\[ \lambda_i=t^{2ij}t^{2i^2+4i} = t^{2i(j+2+i)} = t^{i(2k+4)} = 1^i=1
\]
so we are done.
\end{proof}

\subsection{Preprojective algebras}

Now let $\Ccat=\tl$ and let $\D=S\mModm S$, where $S=\Cnum\times\cdots\times \Cnum$ ($k$ copies) is a basic semisimple $\Cnum$-algebra.  Write $e_i=(0,\ldots,0,1,0,\ldots,0)$ with $1$ in the $i$th entry, so $S$ has $\Cnum$-basis $\{e_1,\ldots,e_k\}$.  The category $\D$ is monoidal by tensoring over $S$. We will consider $\Cnum$-linear monoidal functors $G:\Ccat\to \D$.  Note that $\D\cong \Fun_\Cnum(\M,\M)$, where $\M=S\mMod$, so monoidal functors $G:\Ccat\to \D$ correspond to left module categories $\Ccat\times\M\to\M$: see \cite[Proposition 7.1.3]{egno}.

As $\tl$ is generated by $F_1$, the image of a monoidal functor $G:\Ccat\to \D$ is generated by an $S\da S$-bimodule $X=G(F_1)$.  The category $\D=S\mModm S$ is semisimple and rigid: all simple objects are $1$-dimensional $\Cnum$-vector spaces $S_{ij}$ where $e_iS_{ij}e_j\neq0$, and $S_{ij}^\vee\cong S_{ji}$.  So $X$ is determined up to isomorphism by a quiver $\overline Q$ with vertices $1,\ldots,k$ and an arrow $i\to j$ for each summand $S_{ij}$ of $X$.  As $F_1$ is self-dual we have $X^\vee\cong X$, so $\overline Q$ has the same number of arrows $i\to j$ as $j\to i$.  Therefore $X$ is in fact determined by an (unoriented) graph on the vertices $1,\ldots,k$.

The fusion rule $F_1\otimes F_i\cong F_{i-1}\oplus F_{i+1}$ (see Theorem \ref{thm:tlss}) puts further restrictions on the possible graphs that can appear, and it can be shown (see \cite[Section 6]{ostrik}, \cite[Section 3.3.2]{cooper} and \cite[Chapter I]{ghj}) that the only possibilities are the simply laced Dynkin diagrams of types $A_k$, $D_k$, $E_k$ and the tadpole graph $T_k$.

We now restrict to the $ADE$ cases, excluding the tadpole graphs.  Let $Q$ be a quiver whose underlying graph $\Delta$, where we replace each arrow with an unoriented edge, is an $ADE$ Dynkin diagram.  For each arrow $a:i\to j$ we add an arrow $a^*:j\to i$ to obtain a new quiver $\overline Q$, as above.  Following \cite{gp,dr}, the preprojective algebra $\Pi$ of $\Delta$ is the quotient of the path algebra $\Cnum \overline Q$ by the relations $\sum_{a\in Q_1}aa^*-a^*a$; this is independent of the orientation of $Q$.  The algebra $\Pi$ can also be constructed directly from the Auslander-Reiten theory of $\Cnum Q$, as in \cite{bgl}.

Recall the algebra object $\Sigma$ in $\tl$ from Definition \ref{def:sigma}.  The following result seems to have been known in folklore; a careful statement and proof can be found in \cite[Proposition 5.6.13]{cooper}.
\begin{theorem}\label{thm:tl-pi}
If $G:\tl\to S\mModm S$ is a $\Cnum$-linear monoidal functor with $ADE$ Dynkin graph $\Delta$ then $G(\Sigma)$ is the preprojective algebra of type $\Delta$.
\end{theorem}

As a consequence of Theorem \ref{thm:tl-pi} together with Lemmas \ref{lem:nak-transf} and \ref{lem:sigma-frob} and Theorem \ref{thm:alpha^2}, we get the following:
\begin{corollary}\label{cor:pi-nak2}
Let $\Delta$ be an $ADE$ graph on vertices $1,\ldots,k$ with $k\geq2$.  The preprojective algebra of type $\Delta$, considered as an algebra object in $S\mModm S$, has a Frobenius structure with Nakayama morphism of order $2$.
\end{corollary}

\begin{remark}
Here, we view $\Pi$ as an $S$-algebra with Nakayama morphism as in Definition \ref{def:nak}.  Classically, one views $\Pi$ as a $\Cnum$-algebra and uses the definition of the Nakayama automorphism from \cite[Theorem 1]{na2}.  In this classical setting, Brenner, Butler, and King showed that the Nakayama automorphism has order $2$ \cite[Theorem 4.8]{bbk}.
\end{remark}

\begin{remark}
The fact that the Nakayama automorphism of $\Pi$ has order $2$ can be used to show that the path algebra $\Cnum Q$ of the Dynkin quiver $Q$ is fractionally Calabi-Yau: see \cite[Section 6.2.4]{g-ser}.  This was originally proved by Miyachi and Yekutieli \cite{my}.
\end{remark}

\begin{remark}
The tensor category $\tl$ is related to the representation theory of the Lie algebra $\mathfrak{sl}_2$, as can be seen in its fusion rule.  One can construct algebras from tensor categories with $\mathfrak{sl}_3$ fusion and there is work showing that they have Nakayama automorphisms of order 3: see \cite{ep-nak} together with the 44 pages of calculations in the appendix to the arXiv preprint.  We hope that the methods developed here can be applied in this setting to make calculations more tractable.
\end{remark}

\subsection*{Competing interests}

The authors declare none.


\end{document}